\documentclass[english,12pt,twoside]{article}
\usepackage{amssymb,amsbsy,amsmath,amsfonts,amssymb,amscd,stmaryrd,floatflt}
\usepackage{latexsym,eucal,exscale,epic,eepic,epsfig,amsthm}
\usepackage[all]{xy}
\usepackage{babel}
\usepackage{times}

\newtheorem{theo}{Theorem}[section]

\newtheorem{defi}{Definition}[section]
\newtheorem{lem}{Lemma}[subsection]
\newtheorem{prop}{Proposition}[section]
\newtheorem{rem}{Remark}[section]

\newtheorem{cor}{Corollary}[subsection]

\newcommand{\Sp}{\mathrm Sp}
\newcommand{\Spk}{{\mathrm Sp}^{k}(M)}

\newcommand{\Hom}{\mathrm{Hom}}

\newcommand{\hatI}{\widehat{\mathfrak I}}
\newcommand{\I}{\mathfrak{I}}

\newcommand{\W}{\overline{W}}
\newcommand{\G}{\widetilde{G}}
\renewcommand{\P}{\widetilde{P}}
\newcommand{\PI}{\widetilde{P}_{I}}

\newcommand{\h}{\mathfrak h}

\newcommand{\De}{\Delta}

\newcommand{\ph}{\varphi}
\newcommand{\si}{\sigma}

\newcommand{\la}{\lambda}
\newcommand{\La}{\Lambda}

\begin{document}
\thispagestyle{empty}

\title{Cohomology of Drinfeld symmetric spaces and Harmonic cochains}

\author{Yacine A\"IT AMRANE}

\maketitle

\selectlanguage{english}
\begin{abstract}
Let $K$ be a non-archimedean local field. This paper gives an explicit isomorphism between the dual of the special representation of $GL_{n+1}(K)$ and the space of harmonic cochains defined on the Bruhat-Tits building of  $GL_{n+1}(K)$, in the sense of E. de Shalit \cite{deShalit}. We deduce, applying a work of P. Schneider and U. Stuhler, \cite{Schneider}, that for $K$ of any characteristic, there exists a $GL_{n+1}(K)$-equivariant isomorphism between the cohomology group of the Drinfeld symmetric space and the space of harmonic cochains.
\end{abstract}
{\bf Introduction}\\
Let $K$ be a non-archimedean local field, i.e. a finite extension of ${\mathbb Q}_{p}$ or ${\mathbb F}_{p}((t))$. Let n fixed. Let $\G$ denote $GL_{n+1}(K)$, let $\P$ be its upper triangular Borel subgroup, and let $\overline{S}$ denote the set of fundamental reflexions $s_{i}$, $1\leq i\leq n$, in the linear Weyl group $\W$ of $\G$. Let $\Delta=\{1,\ldots ,n\}$. For each subset $I$ of $\De$, let $\P_{I}$ be the parabolic subgroup of $\G$ generated by $\P$ and the reflexions $s_{i}$, $i\in I$.   

Let $M$ be a commutative ring on which $\G$ acts trivially. For any $I\subseteq \De$, we denote by $C^{^\infty}(\G/\P_{I}, M)$ the space of locally constant functions on $\G/\P_{I}$ with values in $M$. The action of $\G$ on $C^{^\infty}(\G/\P_{I}, M)$ comes from left translations on $\G/\P_{I}$. For any integer $k$, $0\leq k \leq n$, if $J_{k}$ denotes the subset $\{1,\ldots ,n-k\}$ of $\De$, the $k$-special representation of $\G$ is defined to be the $M[\G]$-module :
\vspace{-0.1cm}
$$
\Sp^{k}(M) =\frac{C^{^\infty}(\G/\P_{J_{k}}, M)}{\sum_{j=n-k+1}^{n}C^{^\infty}(\G/\P_{J_{k}\cup\{j\}}, M)}.
$$
In case $k=n$, we get the ordinary Steinberg representation $\Sp^{n}(M)=\textrm{St}^{n}(M)$.

The $n$-dimensional Drinfeld symmetric space over $K$ is the complement $\Omega^{(n+1)}$ in ${\mathbb P}^{n}$ of the union of all the $K$-rational hyperplanes. The group $\G$ acts on $\Omega^{(n+1)}$. 

The symmetric space $\Omega^{(n+1)}$ has been introduced by Drinfeld, \cite{Drinfeld}, who showed that it is endowed with a structure of a rigid analytic variety. In the one dimensional case ($n=1$) when $K$ is of positive characteristic $p>0$, Drinfeld computed the first \'etale cohomology group of $\Omega^{(2)}$ and proved that there are $\G$-isomorphisms :   
\begin{equation}\label{DR}
H_{\textrm{et}}^{^1}(\Omega^{(2)}\otimes_{K}{\bf C}, L) \cong \textrm{Hom}(\textrm{St}^{^{1}}({\mathbb Z}),L) \cong {\frak Harm}^{^{1,1}}({\mathbb Z},L)
\end{equation}
where ${\bf C}$ is the completion of an algebraic closure of $K$, $L$ a finite abelian group whose order is prime to $p$, and  ${\frak Harm}^{^{1,1}}({\mathbb Z},L)$ is the space of $L$-valued harmonic cochains defined on the oriented (or pointed) edges of the Bruhat-Tits tree. (\cite{Drinfeld}, see also \cite{Reversat2}).

Later, in their paper \cite{Schneider}, P. Schneider and U. Stuhler generalized the first isomorphism in (\ref{DR}) to the case of any characteristic of the base field and to any dimension. Indeed, they studied the cohomology groups of $\Omega^{(n+1)}$ for any cohomology theory satisfying certain natural axioms. They proved the existence of a canonical $\G$-equivariant isomorphism, cf. \cite[\S 4, Cor.17]{Schneider}) : 
\begin{equation}\label{SS}
SS :\; H^{\bullet}(\Omega^{(n+1)},{\cal F}) \cong {\rm Hom}_{\mathbb Z}({\rm Sp}^{\bullet}(\mathbb Z),\,L)
\end{equation}
where $\cal F$ is a complex of sheaves on the category of smooth separated rigid analytic varieties over $K$ equipped with a suitable Grothendieck topology, and $L$ is the cohomology of the point $H^{0}({\rm Spec}(K),{\cal F})$. 

If $K$ is of characteristic zero, the isomorphism of Schneider and Stuhler above, applied to rigid De-Rham cohomology, gives a $\G$-isomorphism 
\begin{equation}\label{SSdR}
SS_{dR} :\; H_{dR}^{\bullet}(\Omega^{(n+1)}) \xrightarrow{\simeq}{\rm Hom}_{\mathbb Z}({\rm Sp}^{\bullet}(\mathbb Z),\,K).
\end{equation}

Let $M$ be a commutative ring as above. Let $L$ be an $M$-module on which $\G$ acts linearly. For each $k$, $0\leq k\leq n$, denote by ${\frak Harm}^{k}(M,L)$ the space of $L$-valued harmonic cochains defined over the free $M$-module generated by the pointed $k$-cells of the Bruhat-Tits building associated to $\G$, see def. \ref{cocycles-harmoniques}. In zero characteristic, E. de Shalit, who introduced in \cite{deShalit} the notion of harmonic cochains we use here, proved that there is a $\G$-equivariant isomorphism :
\begin{equation}\label{dSdR}
dS: \; H_{dR}^{\bullet}({\Omega}^{(n+1)}) \xrightarrow{\simeq} {\frak Harm}^{\bullet}({\mathbb Z},K).
\end{equation}
This isomorphism, together with the isomorphism (\ref{SSdR}), gives a $\G$-equivariant isomorphism in characteristic zero : 
\begin{equation}\label{dS}
SS_{dR} \circ dS^{-1}: \; {\frak Harm}^{\bullet}({\mathbb Z},K) \xrightarrow{\simeq} {\rm Hom}_{\mathbb Z}({\rm Sp}^{\bullet}(\mathbb Z),\,K).
\end{equation}
 
In this paper, we shall construct explicitly, in {\bf any characteristic}, the isomorphism (\ref{dS}) above between the harmonic cochain spaces and the $K$-dual spaces of the special representations.
 
The main result in this paper is the following theorem which generalizes also the second Drinfeld's isomorphism in (\ref{DR}) to any dimension $n$ :     

{\bf Theorem \ref{maintheorem} :} Let $K$ be a non-archimedean local field of arbitrary characteristic. Let $M$ and $L$ be as above. Then, for each $k$, $0\leq k\leq n$, there is an explicit $\G$-equivariant isomorphism :
$$
{\frak Harm}^{k}(M,L) \cong {\rm Hom}_{M}({\rm Sp}^{k}(M),\, L).
$$

As a corollary, together with the isomorphism (\ref{SS}) above we obtain the following :

{\bf Corollary :} Let $K$ be a non-archimedean local field of arbitrary characteristic. Let $\cal F$ and $L$ be as in the situation of the isomorphism (\ref{SS}), and as in \cite{Schneider}. For any $k$, $0\leq k\leq n$, we have the following $\G$-equivariant isomorphism : 
$$ 
H^{k}(\Omega^{(n+1)}, {\cal F})\cong {\frak Harm}^{k}({\mathbb Z},L).
$$
In particular, in the case of \'etale cohomology, this isomorphism allows us to express the \'etale cohomology groups of $\Omega^{(n+1)}$ in terms of harmonic cochains which are of combinatorial nature.

Let's summarize the $\G$-isomorphisms we have seen so far by the following commutative diagrams :
$$
\begin{array}{l|r}
\underline{K\;\textrm{of any characteristic}}& \underline{\textrm{car}(K)=0}\;\quad\\
 & \\
\xymatrix{ & \Hom_{\mathbb Z}(\Sp^{k}({\mathbb Z}), L) \\
H^{k}(\Omega^{(n+1)},{\cal F}) \!\!\!\!\!\!\!\!\!\!\!\ar^{SS}[ur] \ar[dr]_{Cor.} & \\
& {\frak Harm}^{k}({\mathbb Z},L) \ar[uu]_{Th. \ref{maintheorem}}}       
& \xymatrix{ \Hom_{\mathbb Z}(\Sp^{k}({\mathbb Z}),\, 
K)  & \\ 
 & \!\!\!\!\!\!\!\!\! H^{k}_{\rm dR}(\Omega^{(n+1)})\;\;\;\;\;\;\;\;\; \ar[ul]_{SS_{dR}} \ar[dl]^{dS}\\
{\frak Harm}^{k}({\mathbb Z},K) \ar[uu]^{SS_{dR}\circ dS^{-1}} & 
}\\
 & \\
 & \\
\textrm{where}\;\; L=H^{0}({\rm Spec}(K),{\cal F}) & 
\textrm{we have}\;\; K=H_{\rm dR}^{0}({\rm Spec}(K)),\;\quad 
\end{array}
$$

The results of this paper were announced without proofs in \cite{Yacine1}. The reader will find more detailled proofs in \cite{Yacine}.  \bigskip\\
Here is the plan of this article. We use the notations introduced above. 

In the first section we give some preliminaries about the Bruhat-Tits building associated to $\G$. For each $I\subseteq \De$, let $B_{I}$ be the standard parahoric subgroup of $\G$ generated by the upper Iwahori subgroup of $\G$ and the fundamental reflexions $s_{i}$, $i\in I$. By the Bruhat decomposition, there is a correspondence between the double classes in $\overline{W}$ and the Bruhat cells in $\G$. By the Iwasawa decomposition and with techniques inspired from Bourbaki \cite{Bourbaki}, we prove that there is a canonical one-to-one correspondence between the double classes in $\overline{W}$ and the Iwasawa cells in $\G$.  By these correspondences and by using decompositions in the Weyl group $\W$ into double classes modulo special subgroups, we deduce decompositions of certain subsets of $\G$ into a disjoint union of Bruhat and Iwasawa cells respectively. 

In the second section, we recall the definition of harmonic cochains given by E. de Shalit. We also recall the relationship, given by  P. Schneider and U. Stuhler, in \cite{Schneider}, between the special representations and the parahoric subgroups. Next, we define, for each $I\subseteq \De$, a subset $C_{I}$ of $\G$ which is a product of standard parahoric subgroups. Finally, by using the Iwasawa decomposition, we prove that the characteristic functions of the open compact subsets $C_{I}\P_{J_{k}}/\P_{J_{k}}$ of $\G/\P_{J_{k}}$, $I\subseteq \De$, viewed in $\Sp^{k}(M)$, have properties that are close to those of harmonic cochains.  

In section 3, we prove the main theorem which gives an explicit isomorphism between duals of the special representations and harmonic cochain spaces. In this isomorphism, the characteristic functions of the subsets $C_{I}\P_{J_{k}}/\P_{J_{k}}$ correspond to the standard cells $\si_{I}$ pointed at the fundamental vertex, stabilized by the standard parahoric subgroups $B_{I}$ (under the action of $\G$ on its Bruhat-Tits building).

\section{Bruhat-Tits building and decompositions in $\G$}
\subsection{Bruhat-Tits building}

From now on, $K$ will be a non-archimedean local field, $O$ its valuation ring, $\pi$ a uniformizing parameter and $\kappa$ the residue field of $K$. We denote by $\G$ the $K$-valued points of the connected reductive linear algebraic group $GL_{n+1}$.

For general properties of buildings, see \cite{Brown} and \cite{Garrett}. An introduction to the Bruhat-Tits building of $\G$ with pointed cells is given  in \cite{deShalit}. \medskip \\    
{\bf The Bruhat-Tits building (pointed cells).} Let $V$ be the $n+1$ vector space $K^{n+1}$. A lattice in $V$ is a free $O$-submodule $\La$ of $V$ of rank $n+1$. The bruhat-Tits building of $\G$ may be described as a simplicial complex $\I$ whose vertices are the dilation classes of lattices $v=[\La]$. Two lattices $\La$ and $\La'$ are in the same class if $\La'=\la\La$ for some $\la \in K^{^*}$. For $k$, $0\leq k\leq n$, a $k$-cell $\si$ in $\I$ is a set of $k+1$ vertices $\{[\La_{0}],[\La_{1}],\ldots ,[\La_{k}]\}$ such that :
\begin{equation}\label{simplexe}
\cdots \supsetneq \La_{0}\supsetneq \La_{1} \supsetneq \cdots \supsetneq \La_{k}\supsetneq \pi\La_{0}\supsetneq \cdots
\end{equation}
Notice that there is an obvious cyclic ordering (mod. $(k+1)$) on the vertices of $\si$.

A pointed $k$-cell of $\I$ is a pair $(\si,v)$ consisting of a $k$-cell $\si$ together with a distinguished vertex $v$ of $\si$. Notice, therefore, that in the case of a pointed cell $(\si,v)$ there is a precise ordering on the vertices. If $v=[\La_{0}]$ we write : 
\begin{equation}\label{simplexepointe}
(\sigma,v) =(\Lambda_{0} \supsetneq \Lambda_{1} \supsetneq \cdots \supsetneq \Lambda_{k} \supsetneq \pi\Lambda_{0}).
\end{equation}

For each $k$, $0\leq k\leq n$, let $\hatI^{k}$ be the set of pointed $k$-cells of $\I$. \medskip\\
{\bf The action of $\G$.} For a fixed basis of the vector space $V$, the action of $\G$ on $V$ is given by the matrix product $ug^{-1}$ where $u\in V$ is considered as a line matrix with respect to the basis of $V$. This action induces an action of $\G$ on the vertex set of the building $\I$ by $g.v=[\La g^{-1}]$. Thus, $\G$ acts on the cells by acting on their vertices. \medskip\\ 
{\bf The type of a pointed cell.} (Cf. \cite[\S\,1.1]{deShalit}.) Let $\si=(\La_{0} \supsetneq \La_{1} \supsetneq \cdots \supsetneq \La_{k} \supsetneq \pi\La_{0})\in {\hatI}^{k}$ be a pointed $k$-cell. The type of  $\si$ is defined as follows :
$$
t(\si)=(d_{1},\ldots ,d_{k+1})
$$
where $d_{i}=\textrm{dim}_{\kappa}\, \La_{i-1}/\La_{i}$ for each $i=1,\ldots ,k+1$ (here, we suppose $\La_{k+1}=\pi\La_{0}$). The type of a pointed $k$-cell is preserved by the action of $\G$. Indeed, the action of $\G$ preserves the dimension of the $\kappa$-vector spaces $\La_{i-1}/\La_{i}$. \medskip\\ 
{\bf The standard cells.} Let $\{u_{1},\ldots ,u_{n+1}\}$ be the fundamental basis of $V$. Consider,  for each $i=0\ldots n$, the vertex $v^{o}_{i}=[\La_{i}^{^o}]$ represented by the lattice :    
$$
\La^{^o}_{i}=\pi Ou_{1}\oplus \cdots \oplus \pi Ou_{i}\oplus Ou_{i+1} \oplus \cdots \oplus Ou_{n+1}.
$$ 
Since the $\La_{i}^{^o},\; 0\leq i\leq n$, satisfy (\ref{simplexe}), we have an $n$-cell $\si_{\emptyset}=\{v^{o}_{0},v^{o}_{1},\ldots ,v^{o}_{n}\}$ called the fundamental chamber of $\I$. 

Now, once and for all, fix $\De=\{1,\ldots ,n\}$. For each $I\subseteq \De$ such that $\De -I=\{i_{1}< \cdots <i_{k}\}$, we have a $k$-cell  
\begin{equation}\label{simplexe-standard-I}
\sigma_{I}=\{v^{o}_{0},v^{o}_{i_{1}},\ldots ,v^{o}_{i_{k}}\}.
\end{equation}
The $\si_{I}$, $I\subseteq \De$, are called the standard cells of the Bruhat-Tits building $\I$. These cells are the faces of the fundamental chamber $\si_{\emptyset}$ having $v^{o}_{0}$ as vertex, called the fundamental vertex of $\I$.    
 
We denote by $\widetilde{T}$ the maximal diagonal torus of $\G$ and by $\widetilde{N}$ its normalizer in $\G$. Since the Weyl group $\W=\widetilde{N}/\widetilde{T}$ of $\G$ with respect to $\widetilde{T}$ is isomorphic to the permutation group ${\cal S}_{n+1}$, then $\W$ is generated by the set $\overline{S}=\{s_{i},\, i\in \De\}$ of the reflexions $s_{i}$ which correspond to the transpositions $(i,i+1)\in {\cal S}_{n+1}$. 
We have the following lemma :  
\begin{lem}\label{wij}
Let $y_{i}$, $0\leq i\leq n$, be the diagonal matrix $y_{i}=\rm{diag}(\overbrace{1,\ldots ,1}^{i\hbox{\scriptsize~times}},\pi,\ldots ,\pi)$ and let $w_{i}=(s_{i}s_{i+1}\cdots s_{n})(s_{i-1}s_{i}\cdots s_{n-1})\cdots (s_{1}s_{2}\cdots s_{n-i+1})\in \W$. We have :
$$
(\si_{\emptyset},v_{i}^{o})=y_{i}w_{i}(\si_{\emptyset},v_{0}^{o}). 
$$
If $(\si,v_{i_{j}}^{o})=(v_{i_{j}}^{o},\ldots ,v_{i_{k}}^{o},v_{i_{0}}^{o},v_{i_{1}}^{o},\ldots ,v_{i_{j-1}}^{o})$ is a face of the pointed chamber $(\si_{\emptyset},v_{i_{j}}^{o})$, where $0\leq i_{0}\lneq i_{1}\lneq \cdots \lneq i_{k}\leq n$ and $0\leq j\leq k$, then 
$$
(\si,v_{i_{j}}^{o})=y_{i_{j}}w_{i_{j}}(\si_{{\widehat I}_{i_{j}}},v_{0}^{o})
$$
where $\De-{\widehat I}_{i_{j}}=\{i_{j+1}-i_{j}<\cdots <i_{k}-i_{j}<n+1+i_{0}-i_{j}<\cdots <n+1+i_{j-1}-i_{j}\}$.
\end{lem}
\proof
The vertices of the fundamental chamber are $v^{o}_{l}=[\Lambda_{l}^{^o}]$, we can easily check that the representants $\Lambda^{^o}_{l}$ of these vertices satisfy : 
$$
\Lambda^{^o}_{l}y_{i}w_{i}= 
\left\{ \begin{array}{lll}
\Lambda^{^o}_{n+1+l-i}   &  \textrm{if} & 0\leq l\leq i-1\\
\Lambda^{^o}_{l-i}\pi & \textrm{if} & i\leq l\leq n,
\end{array}
\right.
$$ 
therefore, by taking into account the way in which $\G$ acts on the vertices of $\I$, it follows that : 
$$
w_{i}^{-1}y_{i}^{-1}v^{o}_{l}= 
\left\{ \begin{array}{lll}
v^{o}_{n+1+l-i}   &  \textrm{if} & 0\leq l\leq i-1\\
v^{o}_{l-i} & \textrm{if} & i\leq l\leq n,
\end{array}
\right.
$$
hence $w_{i}^{-1}y_{i}^{-1}(\si_{\emptyset},v_{i}^{o})=(\si_{\emptyset},v_{0}^{o})$ and, if $(\si,v_{i_{j}}^{o})$ and $\widehat{I}_{i_{j}}$ are as in the lemma, we have $w_{i_{j}}^{-1}y_{i_{j}}^{-1}(\sigma,v^{o}_{i_{j}})=(\sigma_{\widehat{I}_{i_{j}}},v^{o}_{0})$. \qed 

Since the action of $\G$ is transitive on the chambers of $\I$, the lemma above shows that $\G$ acts transitively on the pointed $k$-cells of a given type. Furthermore, if we denote by $t_{I}$ the type of the pointed standard $k$-cell $(\si_{I},v_{0}^{o})$ and by $\hatI^{k,t_{I}}$ the set of all pointed $k$-cells of type $t_{I}$, we have $\hatI^{k} = \coprod_{I\subseteq \De}\hatI^{k,t_{I}}$, where the disjoint union is taken over the subsets $I\subseteq \De$ such that $\De -I$ is of cardinal $k$. Notice, therefore, that for $k$ fixed, there are exactly $C_{n}^{k}$ types of pointed $k$-cells. \medskip \\

\begin{rem}\label{rempointe} For each $I\subseteq \De$, let $B_{I}$ be the pointwise stabilizer in $\G$ of the standard cell $\si_{I}$, or equivalentely the stabilizer of the pointed standard cell $(\si_{I},v_{0}^{o})$. The first assertion of the lemma \ref{wij} shows that, for every $i$, $0\leq i\leq n$, we have  
\begin{equation} \label{pointe}
y_{i}w_{i}B=By_{i}w_{i},
\end{equation}
where $B=B_{\emptyset}$. 
\end{rem}
\subsection{Bruhat and Iwasawa decomposition in $\G$}

\subsubsection{The Bruhat decomposition}
{\bf The parabolic subgroups of $\G$.} Let $\P$ be the upper triangular Borel subgroup of $\G$. A parabolic subgroup of $\G$ is a closed subgroup which contains a Borel subgroup. The subgroups which contain $\P$ are said to be special; these subgroups are completely determined by the subsets $I$ of $\De$. Indeed, if for each $I\subseteq \De$, we let $W_{I}$ be the subgroup of $\W$ generated by the $s_{i},\, i\in I$, it has been shown that the subset 
$$
\P_{I}=\P W_{I}\P \qquad (:=\widetilde{P}\widetilde{N}_{I}\widetilde{P}\quad\textrm{where}\quad  \widetilde{N}_{I}\subseteq \widetilde{N} \;\; \textrm{is such that} \;\;   \widetilde{N}_{I}/\widetilde{T}=W_{I}) 
$$
is a subgroup of $\G$ containing $\P$, and that every subgroup of $\G$ containing $\P$ is a certain $\PI$ for $I\subseteq \De$. Note that $\P=\P_{\emptyset}$.\medskip\\
{\bf The parahoric subgroups of $\G$.} For each $I\subseteq \De$, we denote by $B_{I}^{^\circ}$ the open compact subgroup of $\G(O)$ which is the inverse image of the standard parabolic subgroup $\P_{I}(\kappa)$ of $\G(\kappa)$ by the map ``reduction mod. $\pi$'' : $\G(O)\rightarrow \G(\kappa)$. The parahoric subgroups of $\G$ are  the conjugates in $\G$ of the $B_{I}^{^\circ}$, $I\subseteq \De$. Note that we have $B_{I}=B_{I}^{^\circ}K^{*}$.\label{parahoric-subgroups}\medskip\\
{\bf The Bruhat decomposition.} Let's recall, for each $I\subseteq \De$, the following Bruhat decomposition (cf. \cite[ch. V]{Brown}, \cite{Bourbaki} or \cite{Garrett}) : 
\begin{equation}
B_{I}=BW_{I}B=\coprod_{w\in W_{I}}BwB \qquad \textrm{resp.}\qquad \widetilde{P}_{I}=\widetilde{P}W_{I}\widetilde{P}=\coprod_{w\in W_{I}}\widetilde{P}w\widetilde{P}.
\end{equation}
As a consequence of the Bruhat decomposition, we obtain the following proposition :  
\begin{prop}\label{Bruhatrestriction}
Let $I_{1},I_{2} \subseteq \De$. The map which to $W_{I_{1}}wW_{I_{2}}$ associates $B_{I_{1}}wB_{I_{2}}$ for $w\in \W$  is a one-to-one correspondence :
$$
W_{I_{1}}\backslash{\W}/W_{I_{2}} \xrightarrow{\sim} B_{I_{1}}\backslash \G(O)K^{^*}/B_{I_{2}}.
$$
\end{prop}
\proof Cf. \cite[ch. IV, \S\,2.5, rem. 2]{Bourbaki}.\qed

\subsubsection{The Iwasawa decomposition}

In the following, we shall use the same techniques as in \cite[ch. IV, \S\,2,2]{Bourbaki} and use the generalized Iwasawa decomposition (see for example \cite[th. 17.6]{Garrett}) :
\begin{equation}\label{Iwasawa}
\widetilde{G}=\coprod_{w\in \overline{W}}Bw\widetilde{P}
\end{equation}
to prove the theorem \ref{I1I2} below which gives an analogous result to the proposition \ref{Bruhatrestriction}.

\begin{lem}\label{bourbaki1}
Let $w \in \overline{W}$ and $j\in \Delta$, we have the following inclusions :
\par 1. $w\widetilde{P}s_{j} \subseteq Bw\widetilde{P} \cup Bws_{j}\widetilde{P}$
\par 2. $s_{j}Bw \subseteq Bw\widetilde{P} \cup Bs_{j}w\widetilde{P}$.
\end{lem}
\proof  Indeed, by putting $B'=w^{-1}Bw$ in the first inclusion (resp. $\widetilde{P}'=w\widetilde{P}w^{-1}$ in the second inclusion) we have to show :
$$
\widetilde{P}s_{j} \subseteq B'\widetilde{P} \cup B's_{j}\widetilde{P} \qquad (\textrm{resp.}\; s_{j}B\subseteq B\widetilde{P}' \cup Bs_{j}\widetilde{P}').
$$
The canonical basis of $K^{n+1}$ being $\{u_{1},\ldots ,u_{n+1}\}$, let $\widetilde{G}_{j}$ be the subgroup of  $\widetilde{G}$ consisting of the elements which fix the $u_{i}$ for $i\neq j,j+1$ and which fix the plane spanned by $u_{j}$ and $u_{j+1}$. Put $\widetilde{G}_{j}(O)=\widetilde{G}_{j}\cap \widetilde{G}(O)$. So (cf. \cite[ch.IV, \S 2.2]{Bourbaki}), we have $\widetilde{G}_{j}(k)\widetilde{P}(k)=\widetilde{P}(k)\widetilde{G}_{j}(k)$ for any base field $k$, hence, for $k=K$ (resp. $k=\kappa$) we get $\widetilde{G}_{j}\widetilde{P}=\widetilde{P}\widetilde{G}_{j}$  (resp. $\widetilde{G}_{j}(O)B=B\widetilde{G}_{j}(O)$, by lifting the equality $\widetilde{G}_{j}(\kappa)\widetilde{P}(\kappa)=\widetilde{P}(\kappa)\widetilde{G}_{j}(\kappa)$ to $\widetilde{G}(O)$ and multiplying then by $K^{*}$). Therefore, it's enough to prove : 
$$
\widetilde{G}_{j}\subseteq (B'\cap \widetilde{G}_{j})( \widetilde{P}\cap  \widetilde{G}_{j}) \cup (B'\cap \widetilde{G}_{j})s_{j}( \widetilde{P}\cap  \widetilde{G}_{j})\;
$$
$$
\qquad (\textrm{resp.}\quad
\widetilde{G}_{j}\subseteq (B\cap \widetilde{G}_{j})(\widetilde{P}'\cap \widetilde{G}_{j})\cup (B\cap \widetilde{G}_{j})s_{j}(\widetilde{P}'\cap \widetilde{G}_{j})\,).\qquad\qquad\quad
$$
By identifying $\widetilde{G}_{j}$ with $GL_{2}$, the proof may be completed as in [loc. cit.], except that we use the Iwasawa decomposition instead of the Bruhat decomposition.  \qed

\begin{cor}
\label{bourbaki2}
Let $u_{1},\ldots ,u_{d} \in \overline{S}$ and $w\in \overline{W}$. We have :
\begin{enumerate}
\item $w\widetilde{P}u_{1}\ldots u_{d} \subseteq \displaystyle\bigcup_{(l_{1},\ldots ,l_{p})}Bwu_{l_{1}}\ldots u_{l_{p}}\widetilde{P}$
\item $u_{1}\ldots u_{d}Bw \subseteq \displaystyle\bigcup_{(l_{1},\ldots ,l_{p})}Bu_{l_{1}}\ldots u_{l_{p}}w\widetilde{P}$
\end{enumerate}
where $(l_{1},\ldots ,l_{p})$ runs through the increasing sequenses (including the empty sequence) in $\llbracket 1,d \rrbracket$.
\end{cor}
\proof Induct on $d$ and use the lemma \ref{bourbaki1} above (see also \cite[ch. IV, \S\,2, lem. 1] {Bourbaki}).\qed

\begin{cor}\label{bourbaki3}
Let $I_{1},I_{2}\subseteq \De$. For each $ w\in \W$, we have $B_{I_{1}}w\P_{I_{2}}= BW_{I_{1}}wW_{I_{2}}\P$.
\end{cor}
\proof
Let $I_{1},I_{2}$ and $w$ be as above. Let $w'=u'_{1}\cdots u'_{d_{1}} \in W_{I_{1}}$ and $w''=u''_{1}\cdots u''_{d_{2}}\in W_{I_{2}}$. We have :
$$
Bw'B.Bw\widetilde{P}.\widetilde{P}w''\widetilde{P}=Bu'_{1}\cdots u'_{d_{1}}Bw\widetilde{P}u''_{1}\cdots u''_{d_{2}}\widetilde{P},
$$
therefore, the corollary \ref{bourbaki2} gives $Bw'B.Bw\widetilde{P}.\widetilde{P}w''\widetilde{P}\subseteq BW_{I_{1}}wW_{I_{2}}\widetilde{P}$, and if we take the union as $w'$ and $w''$ run through $W_{I_{1}}$ and $W_{I_{2}}$ respectively, one obtains :
$$
B_{I_{1}}w\widetilde{P}_{I_{2}}\subseteq BW_{I_{1}}wW_{I_{2}}\widetilde{P}.
$$ 
The other inclusion is obvious. \qed

\begin{theo}\label{I1I2}
Let $I_{1},I_{2}\subseteq \De$. The map $\overline{W}\rightarrow B_{I_{1}}\backslash \widetilde{G}/\widetilde{P}_{I_{2}}$ which to $w$ associates $B_{I_{1}}w\widetilde{P}_{I_{2}}$, induces a one-to-one map : 
$$
W_{I_{1}}\backslash \overline{W}/W_{I_{2}} \xrightarrow{\sim} B_{I_{1}}\backslash \widetilde{G}/\widetilde{P}_{I_{2}}.
$$
\end{theo}
\proof The generalized Iwasawa decomposition (\ref{Iwasawa}) shows that the map $w\mapsto Bw\widetilde{P}$ is bijective from $\overline{W}$ on the set $B\backslash \widetilde{G}/\widetilde{P}$, so, by the corollary \ref{bourbaki3}, the surjective map $\overline{W} \xrightarrow{\sim} B\backslash \widetilde{G}/\widetilde{P} \twoheadrightarrow   B_{I_{1}}\backslash \widetilde{G}/\widetilde{P}_{I_{2}}$ induces the following surjective map :
$$
W_{I_{1}}\backslash \overline{W}/W_{I_{2}}   \longrightarrow  B_{I_{1}}\backslash \widetilde{G}/\widetilde{P}_{I_{2}}.
$$
In order to  prove that this map is injective, it is enough to prove the following property:
$$
\textrm{for any}\; w,w'\in \overline{W},\;  B_{I_{1}}w\widetilde{P}_{I_{2}}=B_{I_{1}}w'\widetilde{P}_{I_{2}}, \; \textrm{if and only if,}\; W_{I_{1}}wW_{I_{2}}=W_{I_{1}}w'W_{I_{2}}.
$$ 
Indeed, suppose $B_{I_{1}}w\widetilde{P}_{I_{2}}\cap B_{I_{1}}w'\widetilde{P}_{I_{2}}\neq \emptyset$, so there exist $b\in B^{^\circ}_{I_{1}}$ and $p\in \widetilde{P}_{I_{2}}$ with $bwp=w'$. This implies $p=w^{-1}b^{-1}w' \in \widetilde{G}(O)\cap \widetilde{P}_{I_{2}}\subseteq B_{I_{2}}$ and hence $B_{I_{1}}wB_{I_{2}}\cap B_{I_{1}}w'B_{I_{2}}\neq \emptyset$ which, by Proposition \ref{Bruhatrestriction}, gives $W_{I_{1}}wW_{I_{2}}=W_{I_{1}}w'W_{I_{2}}$.\qed

\begin{rem}\label{remB} Let $I_{1},I_{2}\subseteq \De$. Recall, see \cite[ch. IV, \S\,2.5, prop. 2]{Bourbaki}, that for the Bruhat cells, we have a similar formula to the formula in the corollary \ref{bourbaki3}, that is :
\begin{equation}\label{bruhat2}
B_{I_{1}}wB_{I_{2}}= BW_{I_{1}}wW_{I_{2}}B
\end{equation}

Now, let $I_{1},I_{2}\subseteq \De$ such that for each $i\in I_{1}$, $j\in I_{2}$, we have $|i-j|\geq 2$ (which gives $s_{i}s_{j}=s_{j}s_{i}$). Then, since every element in $W_{I_{1}}$ commutes with every element in $W_{I_{2}}$, we get :
$$
W_{I_{1}\cup I_{2}}=W_{I_{1}}.W_{I_{2}}=W_{I_{2}}.W_{I_{1}}.
$$
The equality (\ref{bruhat2}), for $w=1$, gives then : 
\begin{equation}\label{B}
B_{I_{1}\cup I_{2}}=B_{I_{1}}.B_{I_{2}}=B_{I_{2}}.B_{I_{1}}.
\end{equation}
Notice also that, for each $I\subseteq \De$, one gets (see also \cite[lem. 14 (ii), \S\,4]{Schneider}) :
\begin{equation}\label{BP}
B_{I}\widetilde{P}_{I}=B\widetilde{P}_{I}=B_{I}\widetilde{P}.
\end{equation}
\end{rem}

\subsubsection{Decomposition in the Weyl group $\W$}

For each $r,r'\in \Delta=\{1,\ldots ,n\}$ such that $r \leq r'+1$, we set $w_{r}^{r'}=s_{r}s_{r+1}\ldots s_{r'}$, ($w_{r'+1}^{r'}=1$). By using the Coxeter relations in the Weyl group $\overline{W}$ :
\begin{equation} \label{coxrel'}
\left\{\begin{array}{lll}
s_{l}^{2}=1 & \textrm{for} & l=1,\ldots ,n \\
s_{l_{1}}s_{l_{2}}=s_{l_{2}}s_{l_{1}} & \textrm{for} & 1\leq l_{1} < l_{2}-1 \leq n-1 \\
s_{l}s_{l+1}s_{l}=s_{l+1}s_{l}s_{l+1} & \textrm{for} & l=1, \ldots ,n-1,
\end{array}\right.
\end{equation}
it is easy to show that we have : 
\begin{equation}\label{sw}
s_{l}w_{r}^{r'}=w_{r}^{r'}s_{l-1} \qquad \textrm{for every}\quad r,r'\;\textrm{and}\; l \quad \textrm{such that} \quad r+1 \leq l\leq r'\leq n.
\end{equation}

In all what follows, for any integers $a$ and $b$ such that $a\leq b+1$, we denote by $\llbracket a,b\rrbracket$ the set of all integers $j$ such that $a\leq j\leq b$. It is the empty set in the case $a=b+1$.
  
\begin{prop}\label{wsj2}
Let $k$ be an integer such that $1\leq k\leq n$. Let $w\in \overline{W}$. Then :
\begin{enumerate}
\item For each integer $j$ such that $n-k+2 \leq j \leq n$, we have the decompositions :
$$
BwB_{J_{k} \cup\{j\}}=BwB_{J_{k}} \amalg Bws_{j}B_{J_{k}}\quad \textrm{and} \quad Bw\widetilde{P}_{J_{k} \cup\{j\}}=Bw\widetilde{P}_{J_{k}} \amalg Bws_{j}\widetilde{P}_{J_{k}}.
$$
\item For $j=n-k+1$, we have the following decompositions :
$$
BwB_{{J_{k}} \cup \{n-k+1\}}=\coprod_{r\in \llbracket 1,n-k+2\rrbracket} Bww_{r}^{n-k+1}B_{J_{k}}
$$ 
and 
$$ Bw\widetilde{P}_{{J_{k}} \cup \{n-k+1\}}=\coprod_{r\in \llbracket 1,n-k+2\rrbracket}Bww_{r}^{n-k+1}\widetilde{P}_{J_{k}}.
$$
\end{enumerate}
\end{prop}
\proof We use the proposition \ref{I1I2} and the theorem \ref{Bruhatrestriction} to conclude respectively the first and the second decomposition in 1., and also in 2., from the decomposition into left cosets in $\overline{W}$. \\
1. It is easy to check that, since we have $n-k+2\leq j\leq n$, the sets $W_{J_{k}}$ and $s_{j}W_{J_{k}}$ are the only different left cosets in $W_{J_{k}\cup \{j\}}$ modulo $W_{J_{k}}$.\\
2. To prove that the only left cosets in $W_{J_{k}\cup\{n-k+1\}}$ modulo $W_{J_{k}}$ are the $w_{r}^{n-k+1}W_{J_{k}}$, with $1\leq r\leq n-k+2$, it suffices to prove more generally that for any $a_{1}$ such that $1\leq a_{1}\leq n-k+2$ we have the property :
$$
\begin{array}{c}
{\mathbf A_{\bf 1}} \; :\; \textrm{For each } w\in W_{\llbracket a_{1},n-k+1\rrbracket} \textrm{, there is an } r\in \llbracket a_{1},n-k+2\rrbracket  \textrm{ so that } \\ \\ wW_{J_{k}}=w_{r}^{n-k+1}W_{J_{k}}.
\end{array}
$$
We do this by induction on the length $l(w)$ of $w$ in $W_{\llbracket a_{1},n-k+1\rrbracket}$, where the length is defined with respect to the generators $\overline{S}$. If $l(w)=0$ then $w=1$, so the equality in $\bf A_{1}$ holds trivially with $r=n-k+2$. Now let $w\in  W_{\llbracket a_{1},n-k+1\rrbracket}$ be such that $l(w)=d+1$. Therefore, if $j\in  \llbracket a_{1},n-k+1\rrbracket$ is such that $l(s_{j}w)=d$, by induction, there is an $r\in \llbracket a_{1},n-k+2\rrbracket$ so that : $wW_{J_{k}}=s_{\iota}s_{\iota}wW_{J_{k}}=s_{\iota}w_{r}^{n-k+1}W_{J_{k}}$. To complete the proof of $\bf A_{1}$, we study several cases depending on $\iota$ and $r$ :
\begin{enumerate}
\item[$\bullet$] $\iota \leq r-2$ : the elements $s_{\iota}$ and $w_{r}^{n-k+1}$ commute and $s_{\iota}\in W_{J_{k}}$. Therefore,  
$$
s_{\iota}w_{r}^{n-k+1}W_{J_{k}}=w_{r}^{n-k+1}s_{\iota}W_{J_{k}}=w_{r}^{n-k+1}W_{J_{k}}.
$$
\item[$\bullet$] $\iota =r-1$ (resp. $\iota =r$) : we have :
$$
s_{\iota}w_{r}^{n-k+1}W_{J_{k}}=w_{r-1}^{n-k+1}W_{J_{k}} \quad (\textrm{resp.}\quad s_{\iota}w_{r}^{n-k+1}W_{J_{k}}=w_{r+1}^{n-k+1}W_{J_{k}}).
$$
\item[$\bullet$] $r+1\leq \iota \leq n-k+1$ : by (\ref{sw}) we have $s_{\iota}w_{r}^{n-k+1}=w_{r}^{n-k+1}s_{\iota -1}$, and since $r \leq \iota -1\leq n-k$ then $s_{\iota -1}\in W_{J_{k}}$. Therefore, 
$$
s_{\iota}w_{r}^{n-k+1}W_{J_{k}}=w_{r}^{n-k+1}s_{\iota -1}W_{J_{k}}=w_{r}^{n-k+1}W_{J_{k}}.
$$
\end{enumerate}
To prove that the left cosets $w_{r}^{n-k+1}$, $1\leq r\leq n-k+2$, are different will be done more generally in the proof of the next proposition.
\qed

\begin{prop} \label{generalwsj1'}
Let $k$, $1\leq k\leq n$. Let $a_{1},\ldots, a_{k}$ be such that $1\leq a_{1}\leq \cdots \leq a_{k}$. Assume, furthermore, that $a_{\iota}\leq n-k+\iota +1$ for any $\iota=1,\ldots,k$. We have the decompositions:
$$
B_{\llbracket a_{k},n\rrbracket}\cdots B_{\llbracket a_{1},n-k+1\rrbracket}B_{J_{k}}=\coprod_{(r_{1},\ldots ,r_{k})}\!\! Bw_{r_{k}}^{n}\cdots  w_{r_{1}}^{n-k+1}B_{J_{k}}
$$
and
$$
B_{\llbracket a_{k},n\rrbracket}\cdots B_{\llbracket a_{1},n-k+1\rrbracket}\widetilde{P}_{J_{k}}=\coprod_{(r_{1},\ldots ,r_{k})}\!\! Bw_{r_{k}}^{n}\cdots  w_{r_{1}}^{n-k+1}\widetilde{P}_{J_{k}}
$$
where $(r_{1},\ldots ,r_{k})$ runs through the set $\prod_{\iota=1}^{k}\llbracket a_{\iota},n-k+\iota+1\rrbracket$.
\end{prop}
\proof As above, by Proposition \ref{I1I2} and Theorem \ref{Bruhatrestriction}, it's enough to prove that one has the following decomposition : 
\begin{equation}\label{du}
W_{\llbracket a_{k},n\rrbracket} \cdots W_{\llbracket a_{1},n-k+1\rrbracket} W_{J_{k}}=\coprod_{(r_{1},\ldots ,r_{k})} w_{r_{k}}^{n}\cdots  w_{r_{1}}^{n-k+1}W_{J_{k}}
\end{equation}
with $(r_{1},\ldots ,r_{k})$ running through the set $\prod_{\iota=1}^{k}\llbracket a_{\iota},n-k+\iota+1\rrbracket$.

To prove the equality, it  suffices to prove by induction on $m$ that, for $m=1,\ldots ,k,$ the following holds :
$$
\begin{array}{c}
{\bf A_{m}}: \textrm{for each }w\in W_{\llbracket a_{m},n-k+m\rrbracket} \textrm{ and each } (r_{1},\!\ldots\!,r_{m-1})\in \!\! \displaystyle\prod_{\iota=1}^{m-1}\llbracket a_{\iota},n-k+\iota+1\rrbracket, \\ \;\; \textrm{ there is an } (r'_{1},\ldots, r'_{m})\in \displaystyle\prod_{\iota=1}^{m}\llbracket a_{\iota},n\!-\!k\!+\!\iota \!+\!1\rrbracket \textrm{ so that we have an equality:} \\ \\ w w_{r_{m-1}}^{n-k+m-1}\cdots w_{r_{1}}^{n-k+1}W_{J_{k}}=w_{r'_{m}}^{n-k+m}\cdots w_{r'_{1}}^{n-k+1}W_{J_{k}}.
\end{array}
$$
The proof of $\bf A_{1}$ is above. Assume that $\bf A_{m}$ holds for $m\leq k-1$ and let us show $\bf A_{m+1}$. We have to prove that for any $w\in W_{\llbracket a_{m+1},n-k+m+1\rrbracket}$ and any $(r_{1},\ldots ,r_{m})\in \prod_{\iota=1}^{m}\llbracket a_{\iota},n-k+\iota+1\rrbracket$, there exists $(r'_{1},\ldots ,r'_{m+1})\in \prod_{\iota=1}^{m+1}\llbracket a_{\iota},n-k+\iota+1\rrbracket$ so that :
\begin{equation}\label{rec2}
w w_{r_{m}}^{n-k+m}\cdots w_{r_{1}}^{n-k+1}W_{J_{k}}=w_{r'_{m+1}}^{n-k+m+1}\cdots w_{r'_{1}}^{n-k+1}W_{J_{k}}.
\end{equation}
We prove (\ref{rec2}) by induction on the length $l(w)$ of $w\in W_{\llbracket a_{m+1},n-k+m+1\rrbracket}$. If $l(w)=0$ then $w=1$ and (\ref{rec2}) holds trivially. Assume that (\ref{rec2}) is true when $l(w)=d$. Let $w\in W_{\llbracket a_{m+1},n-k+m+1\rrbracket}$ be such that $l(w)=d+1$. Therefore, if $j\in \llbracket a_{m+1},n-k+m+1\rrbracket$ is such that $l(s_{j}w)=d$, there exists $(r'_{1},\ldots ,r'_{m+1})$ in $\prod_{\iota=1}^{m+1}\llbracket a_{\iota},n-k+\iota+1\rrbracket$, by the induction hypothesis, so that :
$$
s_{j}w w_{r_{m}}^{n-k+m}\cdots w_{r_{1}}^{n-k+1}W_{J_{k}}=w_{r'_{m+1}}^{n-k+m+1}\cdots w_{r'_{1}}^{n-k+1}W_{J_{k}} 
$$
and hence
\begin{equation}\label{m+1}
w w_{r_{m}}^{n-k+m}\cdots w_{r_{1}}^{n-k+1}W_{J_{k}}=s_{j}w_{r'_{m+1}}^{n-k+m+1}\cdots w_{r'_{1}}^{n-k+1}W_{J_{k}}.
\end{equation}
There are several cases depending on $j$ and $r'_{m+1}$ :
\begin{enumerate}
\item[$\bullet$] $a_{m+1}\leq j\leq r'_{m+1}-2 \leq n-k+m$ : we have $s_{j}w_{r'_{m+1}}^{n-k+m+1}=w_{r'_{m+1}}^{n-k+m+1}s_{j}$. Since $s_{j}w_{r'_{m}}^{n-k+m}\in W_{\llbracket a_{m},n-k+m\rrbracket}$, by induction ($\bf A_{m}$), there exists $(r''_{1},\ldots ,r''_{m})$ in $\prod_{\iota=1}^{m} \llbracket a_{\iota},n-k+\iota +1\rrbracket$ such that :
$$
s_{j}w_{r'_{m}}^{n-k+m}\cdots w_{r'_{1}}^{n-k+1}W_{J_{k}}=w_{r''_{m}}^{n-k+m}\cdots w_{r''_{1}}^{n-k+1}W_{J_{k}}.
$$
\item[$\bullet$] $j=r'_{m+1}-1$ : we have $s_{j}w_{r'_{m+1}}^{n-k+m+1}=w_{r'_{m+1}-1}^{n-k+m+1}$.
\item[$\bullet$] $j=r'_{m+1}\leq n-k+m+1$ : we have $s_{j}w_{r'_{m+1}}^{n-k+m+1}=w_{r'_{m+1}+1}^{n-k+m+1}$.
\item[$\bullet$] $a_{m+1}+1\leq r'_{m+1}+1\leq j\leq n-k+m+1$ : we have $s_{j}w_{r'_{m+1}}^{n-k+m+1}=w_{r'_{m+1}}^{n-k+m+1}s_{j-1}$. Since $a_{m}\leq a_{m+1}\leq j-1 \leq n-k+m$, we have $s_{j-1}w_{r'_{m}}^{n-k+m}\in W_{\llbracket a_{m},n-k+m\rrbracket}$ and in the same way as in the first case, one gets :
$$
s_{j-1}w_{r'_{m}}^{n-k+m}\cdots w_{r'_{1}}^{n-k+1}W_{J_{k}}=w_{r''_{m}}^{n-k+m}\cdots w_{r''_{1}}^{n-k+1}W_{J_{k}}.
$$
\end{enumerate}
Thus, together with (\ref{m+1}), there exists $(r''_{1},\ldots ,r''_{m+1})$ in $\prod_{\iota=1}^{m+1}\llbracket a_{\iota},n-k+\iota+1\rrbracket$ such that :
$$
w w_{r_{m}}^{n-k+m}\cdots w_{r_{1}}^{n-k+1}W_{J_{k}}=w_{r''_{m+1}}^{n-k+m+1} w_{r''_{m}}^{n-k+m}\cdots w_{r''_{1}}^{n-k+1}W_{J_{k}}.
$$
This completes the proof of (\ref{rec2}) and also the proof of $\bf A_{m}$, $1\leq m\leq k$.

Let us prove now that the union in (\ref{du}) is a disjoint union. Deny and assume that there are two different elements $(r_{1},\ldots ,r_{k})$ and $(r'_{1},\ldots ,r'_{k})$ in $\prod_{\iota=1}^{k}\llbracket a_{\iota}, n-k+\iota +1\rrbracket$ such that :
$$
w^{n}_{r_{k}}\cdots w^{n-k+1}_{r_{1}}W_{J_{k}}=w^{n}_{r'_{k}}\cdots w^{n-k+1}_{r'_{1}}W_{J_{k}}.
$$
Put $j_{0}=\max \{j,\; 1\leq j\leq k \;|\; r_{j}\neq r'_{j}\}$, without loss of generality we can even assume $r_{j_{0}} > r'_{j_{0}}$. Therefore, since $w^{n-k+j_{0}-1}_{r_{j_{0}-1}}\cdots w^{n-k+1}_{r_{1}}\in W_{\llbracket 1,n-k+j_{0}-1\rrbracket}$, by multiplying the  formula above by $(w^{n}_{r_{k}}\cdots w^{n-k+j_{0}}_{r_{j_{0}}})^{-1}$ on the left and by $W_{\llbracket 1,n-k+j_{0}-1\rrbracket}$ on the right, we get :
$$
W_{\llbracket 1,n-k+j_{0}-1\rrbracket}=(w^{n-k+j_{0}}_{r_{j_{0}}})^{-1}w^{n-k+j_{0}}_{r'_{j_{0}}}W_{\llbracket 1,n-k+j_{0}-1\rrbracket},
$$
hence $(w^{n-k+j_{0}}_{r_{j_{0}}})^{-1}w^{n-k+j_{0}}_{r'_{j_{0}}}\in W_{\llbracket 1,n-k+j_{0}-1\rrbracket}$. As we assumed $r_{j_{0}} > r'_{j_{0}}$, it follows by (\ref{sw}) :
$$ 
(w^{n-k+j_{0}}_{r_{j_{0}}})^{-1}w^{n-k+j_{0}}_{r'_{j_{0}}}=w^{n-k+j_{0}-1}_{r'_{j_{0}}}s_{n-k+j_{0}}(w^{n-k+j_{0}-1}_{r_{j_{0}}-1})^{-1}.
$$
As $w^{n-k+j_{0}-1}_{r'_{j_{0}}},(w^{n-k+j_{0}-1}_{r_{j_{0}}-1})^{-1} \in W_{\llbracket 1,n-k+j_{0}-1\rrbracket}$, this implies that $s_{n-k+j_{0}}$ lies in $W_{\llbracket 1,n-k+j_{0}-1\rrbracket}$, a contradiction.    \qed

\section{Harmonic cochains and special representations}  

Through all this section, we fix $M$ a commutative ring and $L$ an $M$-module. Assume that $\widetilde{G}$ acts trivially on $M$ and that $L$ is endowed with a linear $\widetilde{G}$-action. 

\subsection{Harmonic cochains}

This paragraph concerns some technical lemmas which will be useful to prove the main theorem below (Theorem \ref{maintheorem}). Let us recall the definition of harmonic cochains given by E. de Shalit (\cite[def. 3.1]{deShalit}).
 
\begin{defi}\label{cocycles-harmoniques}
Let $k$ be an integer such that $0\leq k\leq n$. A $k$-harmonic cochain with values in the $M$-module $L$ is a homomorphism $\h\in {\Hom}_{M}(M[\hatI^{k}],\,L)$ which satisfies the following conditions :\\
{\bf (HC1)} If $\si=(v_{0},v_{1},\ldots,v_{k}) \in \hatI^{k}$ and $\si'=(v_{1},\ldots , v_{k},v_{0})$, then 
$$
\h(\si)=(-1)^{k}\h(\si').
$$ 
{\bf (HC2)} Fix a pointed $(k-1)$-cell $\eta\in \hatI^{k-1}$, fix a type $t$ of pointed $k$-cells, and consider the set ${\cal B}(\eta, t)=\{\si \in \hatI^{k};\, \eta<\si\;{\textrm and}\;t(\si)=t\}$. Then 
$$
\sum_{\si\in {\cal B}(\eta,t)}\h(\si)=0. 
$$
{\bf (HC3)} Let $k\geq 1$. Fix $\si=(\La_{0}\supsetneq \La_{1}\supsetneq \La_{k}\supsetneq \pi\La_{0})\in \hatI^{k}$ and fix an index $0\leq j\leq k$. Let ${\cal C}(\si,j)$  be the collection of all $\si'= (\La'_{0}\supseteq \La'_{1}\supsetneq \La'_{k}\supsetneq \pi\La'_{0})\in \hatI^{k}$ for which $\La'_{i}=\La_{i}$ if $i\neq j$, $\La_{j}\supsetneq \La'_{j}$ and $\dim_{\kappa}\La'_{j}/\La_{j+1}=1$. Then 
$$
\h(\si)=\sum_{\si'\in {\cal C}(\si,j)}\h(\si').
$$
{\bf (HC4)} Let $\si =(v_{0},v_{1},\ldots ,v_{k+1})\in \hatI^{k+1}$. Let $\si_{j}=(v_{0},\ldots ,{\hat v}_{j},\ldots ,v_{k+1})\in \hatI^{k}$. Then 
$$
\sum_{j=0}^{k}(-1)^{j}\h(\si_{j})=0.
$$
\end{defi}

For any $k$, $0\leq k\leq n$, we denote by ${\frak Harm}^{k}(M,L)$ the space of $k$-harmonic cochains with values in the $M$-module $L$. In case $k=0$, the condition {\bf (HC4)} shows that
\begin{equation}\label{0-harmonic-cochains}
{\frak Harm}^{0}(M,L)\cong L.
\end{equation}
{\bf The action of $\G$.} The action of $\G$ on ${\frak Harm}^{k}(M,L)$ is induced from its natural action on  ${\Hom}_{M}(M[\hatI^{k}],\,L)$.

To shorten notation, for $I\subseteq \De$, $r'\in \De$ and $r\in I\cup\{r'\}$, set $I_{r}^{r'}=(I \cup \{r'\})-\{r\}$. 
\begin{lem}\label{conditionhc3} 
Let $I\subseteq \De$ such that $\De-I=\{i_{1}<\cdots <i_{k}\}$ and let $j$ such that $1\leq j\leq k$. Let $\h\in {\Hom}_{M}(M[\hatI^{k}],\,L)$ satisfy the condition {\bf (HC3)}. Then :
\begin{enumerate} 
\item If $l$ is such that $j+1\leq l\leq k$, then for each $\si\in {\frak C}(\si_{I\cup\{i_{j}\}},l-1)$ there is $g_{l}^{\si}\in \G$, so that we have the following :
$$
\sum_{\si\in {\cal B}(\si_{I\cup\{i_{j}\}},t_{I})}\!\!\!\h(\si)=\!\!\sum_{\si\in {\cal C}(\si_{I\cup\{i_{j}\}},l-1)}\left(\sum_{\si'\in {\cal B}(\si_{I^{i_{l}}_{i_{l+1}-1}\cup\{i_{j}\}},t_{I^{i_{l}}_{i_{l+1}-1}})}\!\!\!\h(g_{l}^{\si}.\si')\!\right).
$$
\item For $j=l$, there is an integer $m$ such that 
$$
\sum_{\sigma\in {\cal B}(\sigma_{I\cup\{i_{j}\}},t_{I})}\h(\si)=m \left(\sum_{\si'\in {\cal B}(\si_{t_{I^{i_{j}}_{i_{j+1}-1}\cup\{i_{j+1}-1\}}},t_{I^{i_{j}}_{i_{j+1}-1}})}\h(\si')\right).
$$
\end{enumerate}
\end{lem}
\proof First, let $l$ be such that $j\leq l\leq k$. Since $\h$ satisfies the condition {\bf (HC3)}, it follows 
\begin{equation}\label{prov01}
\sum_{\si\in {\cal B}(\si_{I\cup\{i_{j}\}},t_{I})}\h(\si)=\sum_{\si\in {\cal B}(\si_{I\cup\{i_{j}\}},t_{I})}\;\sum_{\si'\in {\cal C}(\si,l)}\h(\si').\smallskip
\end{equation}
Now, let us prove the two assertions of the corollary :\\
1. Assume that $j+1\leq l \leq k$. It is not difficult to show that we have : 
\begin{equation}\label{prov9}
\coprod_{\si\in {\cal B}(\si_{I\cup \{i_{j}\}},t_{I})}{\cal C}(\si,l)=\coprod_{\si\in {\cal C}(\si_{I\cup\{i_{j}\}},l-1)}{\cal B}(\si,t_{I^{i_{l}}_{i_{l+1}-1}}).
\end{equation} 
Combining this with (\ref{prov01}) we get  
\begin{equation}\label{30}
\sum_{\si\in {\cal B}(\si_{I\cup\{i_{j}\}},t_{I})}\h(\sigma)=\sum_{\si\in {\cal C}(\si_{I\cup\{i_{j}\}},l-1)}\;\;\sum_{\si'\in {\cal B}(\si,t_{I^{i_{l}}_{i_{l+1}-1}})}\h(\si').
\end{equation}
The action of $\G$ being transitive on the set of pointed cells of a given type, for each $\si\in {\cal C}(\si_{I\cup\{i_{j}\}},l-1)$ there is $g_{l}^{\si}\in \G$ so that $\si =g_{l}^{\si}(\si_{I_{i_{l+1}-1}^{i_{l}}\cup\{i_{j}\}},v_{0}^{o})$, which implies ${\cal B}(\si,t_{I^{i_{l}}_{i_{l+1}-1}})=g_{l}^{\si}.{\cal B}(\si_{I^{i_{l}}_{i_{l+1}-1}\cup\{i_{j}\}}, t_{I^{i_{l}}_{i_{l+1}-1}})$. Consequently (\ref{30}) can be written as follows 
$$
\sum_{\si\in {\cal B}(\si_{I\cup\{i_{j}\}},t_{I})}\h(\sigma)=\sum_{\si\in {\cal C}(\si_{I\cup\{i_{j}\}},l-1)}\left(\sum_{\si'\in {\cal B}(\si_{I^{i_{l}}_{i_{l+1}-1}\cup\{i_{j}\}},t_{I^{i_{l}}_{i_{l+1}-1}})}\h(g_{l}^{\si}.\si')\right).
$$ 
2. Assume that $l=j$. This assertion being trivial if $i_{j+1}=i_{j}+1$ we can suppose $i_{j+1}-i_{j}\geq 2$. It is not difficult to show that we have the following equality by using the definition of the different sets involved :
\begin{equation}\label{prov56}
\bigcup_{\si\in {\cal B}(\si_{I\cup\{i_{j}\}},t_{I})}{\cal C}(\si ,j)={\cal B}(\si_{I\cup\{i_{j}\}},t_{I_{i_{j+1}-1}^{i_{j}}}). 
\end{equation}
Note that for each $\si'\in {\cal B}(\si_{I\cup\{i_{j}\}},t_{I_{i_{j+1}-1}^{i_{j}}})$, there is exactly $m$ distinct cells $\si\in {\cal B}(\si_{I\cup\{i_{j}\}},t_{I})$ so that $\si'\in {\cal C}(\si, j)$. In total, with (\ref{prov01}), we get 
$$
\sum_{\si\in {\cal B}(\si_{I\cup\{i_{j}\}},t_{I})}\h(\si)=m\left(\sum_{\si'\in {\cal B}(\si_{I\cup\{i_{j}\}},t_{I^{i_{j}}_{i_{j+1}-1}})}\h(\si')\right).
$$
Now, to complete the proof notice that we have the obvious equality $I\cup\{i_{j}\}=I^{i_{j}}_{i_{j+1}-1}\cup\{i_{j+1}-1\}$.  \qed

\begin{lem}\label{prov5}
Let $I\subseteq \De$ with $\De -I = \{i_{1}< \cdots < i_{k}\}$ and let $j$ be an integer such that $1\leq j\leq k$. Then
\begin{enumerate}
\item ${\cal B}(\si_{I\cup\{i_{j}\}},t_{I})=B_{I\cup\{i_{j}\}}.(\si_{I},v_{0}^{o})$.
\item ${\cal C}(\si_{I},j)=B_{I}.(\si_{I_{i_{j+1}-1}^{i_{j}}},v_{0}^{o})$.
\end{enumerate} 
\end{lem}
\proof In both equalities, we prove that the left hand side set is contained in the right hand side set. The inverse inclusions are obvious. \\
1. Let $\si \in {\cal B}(\si_{I\cup\{i_{j}\}},t_{I})$. Since $t(\si)=t_{I}$ and since the action of $\G$ is transitive on the pointed cells of a given type, there exists $b\in \G$ such that $\si=b.(\si_{I},v_{0}^{o})$. Therefore, $(\sigma_{I\cup\{i_{j}\}},v^{o}_{0})$ and $b(\sigma_{I\cup\{i_{j}\}},v^{o}_{0})$ are pointed faces of the same cell $\sigma$. Being also of the same type, we necessarily have $(\sigma_{I\cup\{i_{j}\}},v^{o}_{0})=b.(\sigma_{I\cup\{i_{j}\}},v^{o}_{0})$. Hence $b\in B_{I\cup\{i_{j}\}}$.\\
2. Let $\sigma\in {\cal C}(\sigma_{I},j)$, then, $\sigma=\!(v^{o}_{0},\!\ldots\!,v^{o}_{i_{j-1}},v_{i_{j+1}-1},v^{o}_{i_{j+1}},\!\ldots\!,v^{o}_{i_{k}}\!)$ is the pointed cell obtained from $(\sigma_{I},v^{o}_{0})$ by replacing the vertex $v^{o}_{i_{j}}=[\Lambda^{0}_{i_{j}}]$ by another vertex $v_{i_{j+1}-1}=[\Lambda_{i_{j+1}-1}]$ with $\Lambda_{i_{j+1}-1}\subsetneq \Lambda_{i_{j}}^{0}$ and $\textrm{dim}_{\kappa}(\Lambda_{i_{j+1}-1}/\Lambda_{i_{j+1}}^{0})=1$. We have :
$$
t(\sigma)=t_{I_{i_{j+1}-1}^{i_{j}}}.
$$ 
On the other hand, the pointed $(k+1)$-cell  $\sigma'=(v^{o}_{0},\!\ldots\!,v^{o}_{i_{j-\!1}},v^{o}_{i_{j}},v_{i_{j\!+1}\!-1},v^{o}_{i_{j\!+1}},\!\ldots\!,v^{o}_{i_{k}}\!)$ lies in ${\cal B}(\sigma_{I},\textrm{t}_{I\backslash\{i_{j+1}-1\}})$, thus, by (1) above, there exists $b\in B_{I}$ such that $\sigma'=b.(\sigma_{I\backslash\{i_{j+1}-1\}},v^{o}_{0})$. Acting $b^{-1}$ on the following obvious relation :
$$
(\sigma_{I\cup \{i_{j}\}},v^{o}_{0}) < \sigma < \sigma'=b.(\sigma_{I\backslash\{i_{j+1}-1\}},v^{o}_{0})
$$
we obtain 
$$
(\sigma_{I\cup \{i_{j}\}},v^{o}_{0}) < b^{-1}\sigma < (\sigma_{I\backslash\{i_{j+1}-1\}},v^{o}_{0}).
$$ 
Since $t(b^{-1}\si)=t(\si)=t_{I_{i_{j+1}-1}^{i_{j}}}$, this clearly forces $b^{-1}\si=(\si_{I_{i_{j+1}-1}^{i_{j}}},v^{o}_{0})$. \qed

\subsection{Special representations} 
Let $X$ be a locally compact space. We denote by $C^{^\infty}(X,M)$ (resp. $C_{c}^{^\infty}(X,M)$) the space of locally constant functions on $X$ with values in $M$ (resp. those which, moreover, are compactly supported). Notice that if $X$ is compact then we have :
$$
C_{c}^{^\infty}(X,M)=C^{^\infty}(X,M).
$$   
The spaces $C^{^\infty}(X,M)$ and $C_{c}^{^\infty}(X,M)$ are naturally endowed  with $M$-module structures. Recall, cf. \cite[lemma 4]{Borel-S}, that if $X$ is locally compact, metrizable and totally discontinuous space then : 
\begin{equation}\label{loc.compact}
C_{c}^{^\infty}(X,M)=C_{c}^{^\infty}(X,{\mathbb Z})\otimes_{\mathbb Z}M.
\end{equation}

The group $\G$ is endowed with a topological structure which is induced from the topology of the non-archimedean field $K$. We know that, for any $I\subseteq \De$, the homogeneous space $\G/\P_{I}$ is compact. \medskip\\
{\bf The action of $\G$}. For any $I\subseteq \De$, the action of $\G$ on $C^{^\infty}(\G/\P_{I},M)$ and $C_{c}^{^\infty}(\G/B_{I},M)$ is induced by its action by left translations on respectively $\G/\P_{I}$ and $\G/B_{I}$.   

Let $I\subseteq \De$. For any subset $H$ of $\G$, we denote by $\chi_{H\P_{I}}\in C^{^\infty}(\G/\P_{I},M)$ (resp. $\chi_{HB_{I}}\in C_{c}^{^\infty}(\G/B_{I},M)$) the characteristic function of $H\P_{I}/\P_{I}$ (resp. $HB_{I}/B_{I}$).

\begin{prop} (P. Schneider and U. Stuhler) The $M[\G]$-module $C^{^\infty}(\G/\P_{I},M)$ is generated by the characteristic function $\chi_{B_{I}\P_{I}}$.
\end{prop}
\proof See \cite[\S 4, prop. 8' and cor. 9']{Schneider} and use (\ref{loc.compact}) above. \qed

\begin{rem}
For any $I_{1}\subseteq I_{2}\subseteq \De$,  we have natural commutative diagrams of $M[\G]$-monomorphisms  
$$
\begin{array}{lll}
{} & \!\!\!\!\!\!\!\! \!\!\!\!\!\!\!\!C^{^\infty}(\widetilde{G}/\widetilde{P},\,M)\!\!\!\!\!\!\!\! & {}\\
\qquad\qquad\quad\;\nearrow & {} &  \!\!\!\!\!\!\!\!\nwarrow \\
C^{^\infty}(\widetilde{G}/\widetilde{P}_{I_{2}},M) \!\!\!\!\!\!\!\!\!\!\!\!& \rightarrow & \!\!\!\!\!\!\!\!\!\!C^{^\infty}(\widetilde{G}/\widetilde{P}_{I_{1}},M)  
\end{array}\, {\rm and} \,
\begin{array}{lll}
{} & \!\!\!\!\!\!\!\!\!\!\!\!\!\!\!\!C_{c}^{^\infty}(\G/B,\,M)\!\!\!\!\!\!\!\! & {}\\
\qquad\qquad\quad\;\nearrow & {} & \!\!\!\!\!\!\!\!\nwarrow \\
C_{c}^{^\infty}(\G/B_{I_{2}},M) \!\!\!\!\!\!\!\!\!\!\!\!& \rightarrow & \!\!\!\!\!\!\!\!\!\!C_{c}^{^\infty}(\G/B_{I_{1}},M).
\end{array}
$$
\end{rem}

\begin{defi}Let $k$ be an integer with $0\leq k\leq n$ and let $J_{k}$ be the subset $\llbracket 1,n-k\rrbracket$ of $\De$. A $k$-special representation of $\G$ is the $M[\G]$-module :
$$
\Sp^{k}(M)=\frac{C^{^\infty}(\G/\P_{J_{k}},M)}{\sum_{j=n-k+1}^{n}C^{^\infty}(\G/\P_{J_{k}\cup\{j\}},M)}.
$$
\end{defi}
In case $k=n$, this is the ordinary Steinberg representation. Notice also that, in case $k=0$, this is the trivial representation :
\begin{equation}\label{0-special-representation}
\Sp^{0}(M)\cong M.
\end{equation}
{\bf Link to the parahoric subgroups.} In order to interpret the special representation $\Sp^{k}(M)$ in terms of parahoric subgroups, we recall, following \cite[\S 4]{Schneider}, that we have a surjective map : 
$$
\begin{array}{cccc}
H : &  C_{c}^{^\infty}(\widetilde{G}/B,\, M)  & \longrightarrow & C^{^\infty}(\widetilde{G}/{\widetilde{P}},\, M)
\end{array}
$$
defined by  $H(\varphi)=\displaystyle\sum_{g\in \widetilde{G}/B}\varphi(g)g.\chi_{B\widetilde{P}}$. Recall also that this map induces, for any $I\subseteq \De$, a surjective map :
$$
\begin{array}{cccc}
H_{I} : &  C_{c}^{^\infty}(\widetilde{G}/B_{I},\, M)  & \longrightarrow & C^{^\infty}(\widetilde{G}/{\widetilde{P}_{I}},\, M)
\end{array}
$$
whose kernel is the $M[\widetilde{G}]$-submodule of $C_{c}^{^\infty}(\widetilde{G}/B_{I},\, M)$ generated by the functions $\chi_{By_{i}B_{I}}-\chi_{B_{I}}$, $0 \leq i\leq n$. This leads to the following proposition :

\begin{prop}\label{yeap} For any $k$, $0\leq k\leq n$, we have a canonical isomorphism of $M[\G]$-modules
$$
H_{J_{k}}:\, \frac{C_{c}^{^\infty}(\G/B_{J_{k}},M)}{{\mathfrak R}_{J_{k}}} \cong \Sp^{k}(M),
$$
where ${\mathfrak R}_{J_{k}}$ is the $M[\G]$-submodule of $C^{^\infty}(\G/B_{J_{k}},M)$ generated by the functions $\chi_{B_{J_{k}}s_{j}B_{J_{k}}}+\chi_{B_{J_{k}}}$, $n-k+1\leq j\leq n$, and the functions $\chi_{By_{i}B_{J}}\!-\chi_{B_{J}}$, $0\leq i\leq n$.  \qed
\end{prop}

\subsection{Harmonic cochains and special representations}
\subsubsection{Definition of new sets $C_{I}$ of $\G$.}
For each $I\subseteq \De$, for each $r'_{1},\ldots ,r'_{m}\in \De$ and each $r_{1},\ldots ,r_{m}\in I\cup \{r'_{1},\ldots ,r'_{m}\}$, we set :
$$
I^{r'_{1},\ldots ,r'_{m}}_{r_{1},\ldots ,r_{m}}=(I\cup \{r'_{1},\ldots ,r'_{m}\})-\{r_{1},\ldots ,r_{m}\}.
$$

\par Let us fix an integer $k$ such that $1\leq k\leq n$ and set $J_{k}$ the subset $\llbracket 1,n-k\rrbracket$ of $\De$. 

\par Let $I\subseteq \De$ be such that $\Delta - I=\{i_{1}< \cdots < i_{k}\}$. For every $m=1,\ldots ,k$, we necessarily have $i_{m}\leq n-k+m$, therefore the integers $i_{1},\ldots ,i_{m}$ lie in the subset $J_{k}\cup\{n-k+1,\ldots,n-k+m\}=\llbracket 1,n-k+m\rrbracket$ of $\De$. Hence,  for any $m=1,\ldots ,k$, if we put $i_{0}=0$, one can see easily that we have :
\begin{equation}\label{JI}
{J_{k}}^{n-k+1,\ldots ,n-k+m}_{\;\;\;\;\; i_{1}\;\;\;\; ,\ldots ,\;\;\;\; i_{m}}=\left(\coprod_{\iota =1}^{m}\llbracket i_{\iota -1}+1,i_{\iota}-1\rrbracket\right) \;\amalg \llbracket i_{m}+1,n-k+m\rrbracket.
\end{equation}
Moreover, if we put $m=k$ in this formula (\ref{JI}), we can see that ${J_{k}}^{n-k+1,\ldots,n}_{\;\;\;\;\; i_{1}\;\;\;\; ,\ldots ,i_{k}}=I$. 

\par Now, if $I\subseteq \Delta$ is such that $\Delta-I=\{i_{1}<\cdots <i_{k}\}$, we write :
\begin{equation}\label{defci}
C^{^\circ}_{I}=B^{^\circ}_{{J_{k}}^{n-k+1,\ldots,n}_{\;\;\;\;\; i_{1}\;\;\, ,\ldots ,i_{k}}}\!\!\! \cdots B^{^\circ}_{{J_{k}}^{n-k+1}_{\;\;\;\;\; i_{1}}}B^{^\circ}_{J_{k}}      \quad       \textrm{and}     \quad      
C_{I}=B_{{J_{k}}^{n-k+1,\ldots,n}_{\;\;\;\;\; i_{1}\;\;\, ,\ldots ,i_{k}}} \!\!\! \cdots B_{{J_{k}}^{n-k+1}_{\;\;\;\;\; i_{1}}}B_{J_{k}}.
\end{equation}
The set $C^{^\circ}_{I}$ is compact open in $\G$ and we clearly have $C_{I}=C^{^\circ}_{I}K^{*}$, see \S \ref{parahoric-subgroups}. Hence, the set $C_{I}\widetilde{P}_{J_{k}}/\widetilde{P}_{J_{k}}=C^{^\circ}_{I}\widetilde{P}_{J_{k}}/\widetilde{P}_{J_{k}}$ is compact open in the homogeneous space $\widetilde{G}/\widetilde{P}_{J_{k}}$.

\begin{theo}\label{decomp.CIPJ} If for each $I\subseteq \De$ such that $\De-I=\{i_{1}<\cdots < i_{k}\}$, we define ${\mathfrak C}_{I}=\prod_{\iota=1}^{k} \llbracket i_{\iota}+1, n-k+\iota +1\rrbracket$, then we have the following decompositions :
$$
\chi_{C_{I}}=\sum_{\underline{r}\in {\mathfrak C}_{I}}\chi_{Bw_{r_{k}}^{n}\cdots w_{r_{1}}^{n-k+1}B_{J_{k}}} \quad \textrm{and} \quad \chi_{C_{I}\P_{J_{k}}}=\sum_{\underline{r}\in {\mathfrak C}_{I}}\chi_{Bw_{r_{k}}^{n}\cdots w_{r_{1}}^{n-k+1}\P_{J_{k}}}, 
$$
where $\underline{r}$ denote the $k$-tuple $(r_{1},\ldots,r_{k})$. 
\end{theo}
\proof  For any $m=1,\ldots, k$, the expression (\ref{JI}) above shows that ${J_{k}}^{n-k+1,\ldots ,n-k+m}_{\;\;\;\;\; i_{1}\;\;\;\; ,\ldots ,\;\;\;\; i_{m}}$ decomposes as a union of intervals which satisfy pairwise the hypothesis of the assertion (\ref{B}) given in Remark \ref{remB}. Thus, an easy induction on $m$ by using this assertion (\ref{B}) proves that we have : 
\begin{equation}\label{bruhatprep1}
C_{I} = B_{\llbracket i_{k}+1,n\rrbracket} \cdots B_{\llbracket i_{1}+1,n-k+1\rrbracket} B_{J_{k}}.
\end{equation}
Next, since we have $B_{J_{k}}\P_{J_{k}}=B\P_{J_{k}}$ which is given by the assertion (\ref{BP}) of the same remark, we conclude from (\ref{bruhatprep1}) that we also have :
\begin{equation}\label{iwasawaprep2}
C_{I}\P_{J_{k}}=B_{\llbracket i_{k}+1,n\rrbracket} \cdots B_{\llbracket i_{1}+1,n-k+1\rrbracket} \P_{J_{k}}.
\end{equation} 
Finally, by Proposition \ref{generalwsj1'}, we deduce from (\ref{bruhatprep1}) and (\ref{iwasawaprep2}) the following respective decompositions 
$$
C_{I}=\coprod_{\underline{r}\in {\mathfrak C}_{I}}Bw_{r_{k}}^{n}\cdots w_{r_{1}}^{n-k+1}B_{J_{k}} \quad \textrm{and} \quad C_{I}\P_{J_{k}}=\coprod_{\underline{r}\in {\mathfrak C}_{I}}Bw_{r_{k}}^{n}\cdots w_{r_{1}}^{n-k+1}\P_{J_{k}}, 
$$
and the theorem follows. \qed

For convenience, we will call ${\mathfrak C}_{I}$ the index set associated to the decomposition of $C_{I}$.

\subsubsection{About vanishing in $\Sp^{k}(M)$.}
The propositions below give a method which allows us to know whether certain elements vanish in $\Sp^{k}(M)$.
\begin{prop}\label{astuce0}
Let $w\in \overline{W}$ and $w' \in W_{\llbracket n-k+2,n\rrbracket}$ (In case $k=n$, consider $w'\in \overline{W}$). We have :
$$
\chi_{BwB_{J_{k}}}-(-1)^{l(w')}\chi_{Bww'B_{J_{k}}}\in \sum_{j=n-k+1}^{n}C_{c}^{^\infty}(\widetilde{G}/B_{J_{k}\cup\{j\}},\, M)
$$
and 
$$
\chi_{Bw\widetilde{P}_{J_{k}}}-(-1)^{l(w')}\chi_{Bww'\widetilde{P}_{J_{k}}}\in \sum_{j=n-k+1}^{n}C^{^\infty}(\widetilde{G}/\widetilde{P}_{J_{k}\cup\{j\}},\, M).
$$ 
\end{prop}
\proof We prove the second assertion. Let $u_{1},\ldots ,u_{d} \in \overline{S}$ such that $w'=u_{1}\cdots u_{d}$ is a reduced expression ($d=l(w')$). We have : 
$$
\chi_{Bw\widetilde{P}_{J_{k}}}-(-1)^{d} \,\chi_{Bww'\widetilde{P}_{J_{k}}}=\sum_{r=1}^{d}(-1)^{r-1}(\chi_{Bwu_{1}\cdots u_{r-1}\widetilde{P}_{J_{k}}}+\chi_{Bwu_{1}\cdots u_{r}\widetilde{P}_{J_{k}}}).
$$
The expression $w'=u_{1}\cdots u_{d}$ being reduced and since $w'\in W_{\llbracket n-k+2,n\rrbracket}$, we deduce that for each $r=1,\ldots ,d$, there is an integer $j$ such that $n-k+2\leq j\leq n$ (in case $k=n$, $j$ is such that $1\leq j\leq n$) and $u_{r}=s_{j}$. Thus, by Proposition \ref{wsj2}, we have : 
$$
\chi_{Bwu_{1}\cdots u_{r-1}\widetilde{P}_{J_{k}}}+\chi_{Bwu_{1}\cdots u_{r}\widetilde{P}_{J_{k}}}=\chi_{Bwu_{1}\cdots u_{r-1}\widetilde{P}_{J_{k}\cup\{j\}}}\in C^{^\infty}(\widetilde{G}/\widetilde{P}_{J_{k}\cup\{j\}},\, M).
$$
The proof of the first assertion is similar. \qed

\begin{prop}\label{astuce}
Let $w,w' \in \overline{W}$ and let $a,b$ be two integers such that  $1\leq a \leq b \leq n$. Suppose furthermore that we have $s_{b}w'=w's_{b'}$, where $b'$ is an integer such that $n-k+2 \leq b'\leq n$. Then :
$$
\sum_{(r_{1},r_{2})}\chi_{Bww_{r_{2}}^{b}w_{r_{1}}^{b-1}w'B_{J_{k}}}=\sum_{l=0}^{b-a} \sum_{r=a}^{b-l}\chi_{Bww_{r}^{b}w_{b-l}^{b-1}w'B_{J_{k}\cup\{b'\}}}\; \in C_{c}^{^\infty}(\widetilde{G}/B_{J_{k}\cup\{b'\}},M),
$$
and also
$$
\sum_{(r_{1},r_{2})}\chi_{Bww_{r_{2}}^{b}w_{r_{1}}^{b-1}w'\widetilde{P}_{J_{k}}}=\sum_{l=0}^{b-a} \sum_{r=a}^{b-l}\chi_{Bww_{r}^{b}w_{b-l}^{b-1}w'\widetilde{P}_{J_{k}\cup\{b'\}}}\; \in C^{^\infty}(\widetilde{G}/\widetilde{P}_{J_{k}\cup\{b'\}},M)
$$
where the pair $(r_{1},r_{2})$ runs through the set $\llbracket a,b\rrbracket \times \llbracket a,b+1\rrbracket$.
\end{prop}
\proof First, if we put $b=a+m$, by induction on $m\geq 0$ we prove easily that the set $\llbracket a,b\rrbracket \times \llbracket a,b+1\rrbracket$ decomposes into a disjoint union as follows : 
$$
\llbracket a,b\rrbracket \times \llbracket a,b+1\rrbracket =\coprod_{l=0}^{b-a}\big(\llbracket a,b-l\rrbracket \times \{b-l+1\} \amalg \{b-l\}\times \llbracket a,b-l\rrbracket \big).
$$
Now, let us prove the second formula of the proposition. Since $(r_{1},r_{2})$ runs through $\llbracket a,b\rrbracket \times \llbracket a,b+1\rrbracket$, from the decomposition above we get :
$$
\sum_{(r_{1},r_{2})}\chi_{Bww_{r_{2}}^{b}w_{r_{1}}^{b-1}w'\widetilde{P}_{J_{k}}}=\sum_{l=0}^{b-a}(\sum_{r_{1}=a}^{b-l}\chi_{Bww_{b-l+1}^{b}w_{r_{1}}^{b-1}w'\widetilde{P}_{J_{k}}}+\sum_{r_{2}=a}^{b-l}\chi_{Bww_{r_{2}}^{b}w_{b-l}^{b-1}w'\widetilde{P}_{J_{k}}}).
$$
On the other hand, by using the formula (\ref{sw}) together with the hypothesis of the proposition, we obtain : 
$$
w_{b-l+1}^{b}w_{r_{1}}^{b-1}w'=w_{b-l+1}^{b}w_{r_{1}}^{b}s_{b}w'=w_{r_{1}}^{b}w_{b-l}^{b-1}s_{b}w'=w_{r_{1}}^{b}w_{b-l}^{b-1}w's_{b'}.
$$
By replacing in the first sum of the right hand side of the equality above, we obtain :
$$
\sum_{(r_{1},r_{2})}\chi_{Bww_{r_{2}}^{b}w_{r_{1}}^{b-1}w'\widetilde{P}_{J_{k}}}=\sum_{l=0}^{b-a}\sum_{r=a}^{b-l}(\chi_{Bww_{r}^{b}w_{b-l}^{b-1}w's_{b'}\widetilde{P}_{J_{k}}}+\chi_{Bww_{r}^{b}w_{b-l}^{b-1}w'\widetilde{P}_{J_{k}}}).
$$
Finally, since $n-k+2\leq b'\leq n$, by Proposition \ref{wsj2}, we have 
$$ \chi_{Bww_{r}^{b}w_{b-l}^{b-1}w's_{b'}\widetilde{P}_{J_{k}}}+\chi_{Bww_{r}^{b}w_{b-l}^{b-1}w'\widetilde{P}_{J_{k}}}=\chi_{Bww_{r}^{b}w_{b-l}^{b-1}w'\widetilde{P}_{J_{k}\cup\{b'\}}}.
$$ 
The proof of the first formula is similar. \qed \medskip\\

\subsubsection{Harmonicity in $\bf\Sp^{k}(M)$} The following proposition and its corollary below show that the characteristic functions $\chi_{C_{I}}$ and $\chi_{C_{I}\P_{J_{k}}}$ have properties that are somehow similar to those of harmonic cochains.  First, we will need the following technical lemma :

\begin{lem}\label{indicetranslation2} Let $i$ be an integer such that $0\leq i\leq n$. Let $r_{1},\cdots ,r_{k}$ be integers such that for each $j=1,\cdots, k$, we have  $1\leq r_{j}\leq n-k+j+1-i$. Then
\begin{equation}\label{indicetranslation0}
w_{i}w_{r_{k}}^{n}\cdots w_{r_{1}}^{n-k+1}=w_{r_{k}+i}^{n}\cdots w_{r_{1}+i}^{n-k+1}w_{i}^{n-k}\cdots w_{1}^{n-i+1-k},
\end{equation}
where $w_{i}=w_{i}^{n}w_{i-1}^{n-1}\cdots w_{1}^{n-i+1}$ (cf. lemma \ref{wij}).
\end{lem}
\proof
Let's first prove that, if $a,a',b,b'$ are integers satisfying $1\leq a\leq a'\leq b\leq b'\leq n$, we have the equality :
\begin{equation}\label{indicetranslation3}
w_{a}^{b}w_{a'}^{b'}=w_{a'+1}^{b'}w_{a}^{b-1}.
\end{equation}
Indeed, we can write $w_{a}^{b}=w_{a}^{a'}w_{a'+1}^{b}$. On the other hand, by using the formula (\ref{sw}), we get $w_{a'+1}^{b}w_{a'}^{b'}=w_{a'}^{b'}w_{a'}^{b-1}$. Thus, the left hand side of the equality (\ref{indicetranslation3}) can be written as follows :
\begin{equation}\label{indicetranslation4}
w_{a}^{b}w_{a'}^{b'}=w_{a}^{a'}w_{a'}^{b'}w_{a'}^{b-1}.
\end{equation}
Since we have $w_{a}^{a'}w_{a'}^{b'}=w_{a}^{a'-1}w_{a'+1}^{b'}=w_{a'+1}^{b'}w_{a}^{a'-1}$ and $w_{a}^{a'-1}w_{a'}^{b-1}=w_{a}^{b-1}$, it follows that  the right hand sides of the equalities (\ref{indicetranslation3}) and (\ref{indicetranslation4}) are equal.

Now, in order to establish the formula (\ref{indicetranslation0}) in the lemma, apply the identity (\ref{indicetranslation3}) with $w_{a}^{b}=w_{1}^{n-i+1}$ and $w_{a'}^{b'}=w_{r_{j}}^{n-k+j}$ for each $j=k,\ldots ,1$, and in that order. Next, proceed in the same way with $w_{a}^{b}=w_{\iota}^{n-i+\iota}$ for $\iota =2, \ldots ,i$ respectively. \qed

\begin{prop}\label{astuce2}
Let $I\subseteq \De$ with $\De - I =\{i_{1}< \cdots < i_{k}\}$. Let $i_{k+1}$ be such that $i_{k}\lneq i_{k+1}\leq n+1$ and let $\widehat{I}_{1}\subseteq \De$ be such that 
$$
\Delta - \widehat{I}_{1}=\{i_{2}\!-\!i_{1}<\cdots < i_{k}\!-\!i_{1}<i_{k+1}\!-\!i_{1}\}.
$$
In $\Sp^{k}(M)$, we have the equalities (the first equality is seen in $\Sp^{k}(M)$ through the isomorphism $H_{J_{k}}$ given by Proposition \ref{yeap}) : 
\begin{enumerate}
\item $\chi_{C_{I}}=\displaystyle\sum_{\underline{r}\in \mathfrak{C}_{I}^{0}}\chi_{Bw_{r_{k}}^{n}\cdots w_{r_{1}}^{n-k+1}B_{J_{k}}}+\displaystyle\sum_{t=1}^{k-1}(-1)^{k-t-1}\!\!\!\!\!\!\sum_{\underline{r}\in {\mathfrak C}_{I}^{t,k-t-1}}\!\!\!\chi_{Bw_{r_{k}}^{n}\cdots w_{r_{1}}^{n-k+1}B_{J_{k}}} +\displaystyle\sum_{\underline{r}\in \mathfrak{C}_{I}^{k}}\chi_{Bw_{r_{k}}^{n}\cdots w_{r_{1}}^{n-k+1}B_{J_{k}}}$.
\item $\chi_{C_{I}\widetilde{P}_{J_{k}}}=\displaystyle\sum_{\underline{r}\in \mathfrak{C}_{I}^{0}}\chi_{Bw_{r_{k}}^{n}\cdots w_{r_{1}}^{n-k+1}\widetilde{P}_{J_{k}}}+\displaystyle\sum_{t=1}^{k-1}(-1)^{k-t-1}\!\!\!\!\!\!\sum_{\underline{r}\in {\mathfrak C}_{I}^{t,k-t-1}}\!\!\!\chi_{Bw_{r_{k}}^{n}\cdots w_{r_{1}}^{n-k+1}\widetilde{P}_{J_{k}}} +\displaystyle\sum_{\underline{r}\in \mathfrak{C}_{I}^{k}}\chi_{Bw_{r_{k}}^{n}\cdots w_{r_{1}}^{n-k+1}\widetilde{P}_{J_{k}}}$.
\item
$
y_{i_{1}}w_{i_{1}}\chi_{C_{\widehat{I}_{1}}\widetilde{P}_{J_{k}}}=(-1)^{k}\displaystyle\sum_{\underline{r}\in \mathfrak{C}_{I}^{0}}\chi_{Bw_{r_{k}}^{n}\cdots w_{r_{1}}^{n-k+1}\widetilde{P}_{J_{k}}}
$
\end{enumerate}
where we have set $\;{\mathfrak C}_{I}^{0}=\displaystyle\prod_{\iota =1}^{k}\llbracket i_{\iota}+1,i_{\iota +1}\rrbracket\;$, $\;{\mathfrak C}_{I}^{k}=\displaystyle\left(\prod_{\iota =1}^{k-1}\llbracket i_{\iota}+1,n-k+\iota +1\rrbracket \right) \times \llbracket i_{k+1}+1,n+1\rrbracket$, and for each $t$ such that $1\leq t\leq k-1$ we have set
$$
{\frak C}_{I}^{t,k-t-1}\!\!=\!\!\left(\prod_{\iota =1}^{t-1}\llbracket i_{\iota}\! +\! 1,n\! -\! k\! +\!\iota\! +\!1\rrbracket\right) \!\times \!\left(\prod_{\iota =t}^{k-1}\llbracket i_{\iota +1}\! +\! 1,n\! -\! k\! +\! \iota\! +\! 1\rrbracket\right) \times \llbracket i_{k}\! +\! 1,i_{k+1}\rrbracket.
$$ 
\end{prop}

\proof The proofs of (1) and (2) are similar, we will prove (1) and (3).\\
1. In fact, it's enough to prove that the equality holds in $C_{c}^{^\infty}(\G/B_{J_{k}},M)$ modulo $\sum_{j=n-k+1}^{n}C_{c}^{^\infty}(\G/B_{J_{k}\cup\{j\}},M)$.  Recall that for any $I\subseteq \Delta$, to the subset  $C_{I}$ of $\widetilde{G}$ we have associated an ``index set'' ${\frak C}_{I}$ so that we have the following decomposition (cf. Theorem \ref{decomp.CIPJ}) :
$$
\chi_{C_{I}}=\sum_{\underline{r}\in {\frak C}_{I}}\chi_{Bw_{r_{k}}^{n}\cdots w_{r_{1}}^{n-k+1}B_{J_{k}}}.
$$
On the other hand, if ${\frak C}^{0}_{I}$ and ${\frak C}^{k}_{I}$ are as in the statement of this proposition and if for each $t=1,\ldots,k-1$, we put  
$$
{\mathfrak C}_{I}^{t}=\left(\prod_{\iota =1}^{t-1}\llbracket i_{\iota}+1,n-k+\iota +1\rrbracket\right) \times \llbracket i_{t+1}+1,n-k+t+1\rrbracket \times \left(\prod_{\iota =t+1}^{k}\llbracket i_{\iota}+1,i_{\iota +1}\rrbracket\right),
$$
then it is not difficult to show that ${\frak C}_{I}$ is a disjoint union of the  ${\frak C}_{I}^{t}$, $0\leq t\leq k$. Hence :
\begin{equation}\label{haja}
\chi_{C_{I}}=\sum_{t=0}^{k}\sum_{\underline{r}\in {\frak C}_{I}^{t}}\chi_{Bw_{r_{k}}^{n}\cdots w_{r_{1}}^{n-k+1}B_{J_{k}}}.
\end{equation}
Now, let $t$ be such that $1\leq t \leq k-1$. For each $t'=0,\ldots , k-t-2$, put : 
$$ 
\begin{array}{c}
{\frak C}_{I}^{t,t'}=\left(\displaystyle\prod_{\iota =1}^{t-1}\llbracket i_{\iota}+1,n-k+\iota +1\rrbracket\right) \times \left(\displaystyle\prod_{\iota =t}^{t+t'}\llbracket i_{\iota +1}+1,n-k+\iota +1\rrbracket\right) \\\ 
\qquad\times\llbracket i_{t+t'+1}+1,n-k+t+t'+2\rrbracket \times \left(\displaystyle\prod_{\iota =t+t'+2}^{k}\llbracket i_{\iota}+1,i_{\iota +1}\rrbracket\right), \end{array}
$$
and for $t'=n-k-1$, let ${\frak C}_{I}^{t,n-k-1}$ as in the statement of the proposition. Also, for each $t'=0,\ldots ,k-t-1$, put : 
$$
{\frak D}_{I}^{t,t'}\!\!\!=\!\!\left(\displaystyle\prod_{\iota =1}^{t-1}\llbracket i_{\iota}\! +\! 1,n\! -\! k\! +\!\iota \! +\! 1\rrbracket\!\right) \! \times \! \left(\!\displaystyle\prod_{\iota =t}^{t+t'}\llbracket i_{\iota +1}\! +\! 1,n\! -\! k\! +\!\iota \! +\! 1\rrbracket\!\right) \!
\times \! \left(\displaystyle\prod_{\iota \!=t+t'\!+\!1}^{k}\!\!\!\!\!\!\llbracket i_{\iota}\! +\! 1,i_{\iota +1}\rrbracket\!\right),
$$
and for $t'=k-t$, let ${\frak D}_{I}^{t,k-t}=\emptyset$. Notice that the interval which corresponds to $\iota=t+t'+1$ in ${\frak C}_{I}^{t,t'}$ can be decomposed as follows :
$$
\llbracket i_{t+t'+1}+1, n-k+t+t'+2\rrbracket =\llbracket i_{t+t'+1}+1,i_{t+t'+2}\rrbracket  \amalg \llbracket i_{t+t'+2}+1,n-k+t+t'+2\rrbracket, 
$$ 
therefore, for each $t'$, $0\leq t'\leq k-t-1$, we have ${\frak C}_{I}^{t,t'}={\frak D}_{I}^{t,t'}\amalg {\frak D}_{I}^{t,t'+1}$. Consider the alternating sum over $t'$ as follows :
$$
\sum_{t'=0}^{k-t-1}(-1)^{t'}\!\!\sum_{\underline{r}\in {\frak C}_{I}^{t,t'}}\chi_{Bw_{r_{k}}^{n}\cdots w_{r_{1}}^{n-k+1}B_{J_{k}}}=\sum_{t'=0}^{k-t-1}(-1)^{t'}\!\!\!\!\sum_{\underline{r}\in {\frak D}_{I}^{t,t'}\amalg {\frak D}_{I}^{t,t'+1}}\chi_{Bw_{r_{k}}^{n}\cdots w_{r_{1}}^{n-k+1}B_{J_{k}}}.
$$
In the right hand side of this equality, we see that all the sums over the ${\frak D}_{I}^{t,t'}$, $t'=1,\ldots ,k-t-1,$ cancel each other. What remains is the sum over ${\frak D}_{I}^{t,0}={\frak C}_{I}^{t}$ and the sum over ${\frak D}_{I}^{t,k-t}=\emptyset$. Therefore we get :
\begin{equation}\label{frmle2}
\displaystyle\sum_{\underline{r}\in {\frak C}_{I}^{t}} \chi_{Bw_{r_{k}}^{n}\cdots w_{r_{1}}^{n-k+1}B_{J_{k}}}=\sum_{t'=0}^{k-t-1}(-1)^{t'}\sum_{\underline{r}\in {\frak C}_{I}^{t,t'}}\chi_{Bw_{r_{k}}^{n}\cdots w_{r_{1}}^{n-k+1}B_{J_{k}}}.
\end{equation}
Notice that in ${\frak C}_{I}^{t,t'}$, for each $t'=0,\ldots k-t-2$, the two intervals which correspond to the indeces $t+t'$ and $t+t'+1$ are of the form of those in Proposition \ref{astuce} with $a=i_{t+t'+1}+1$ and $b=n-k+t+t'+1$, and if we put $w'=w_{r_{t+t'-1}}^{n-k+t+t'-1}\cdots w_{r_{1}}^{n-k+1}$, it is clear that we have $s_{b}w'=w's_{b}$. Therefore, in $\Sp^{k}(M)$, we have (for each $t=1,\ldots ,k-1$) :
$$
\displaystyle\sum_{\underline{r}\in {\frak C}_{I}^{t}} \chi_{Bw_{r_{k}}^{n}\cdots w_{r_{1}}^{n-k+1}B_{J_{k}}}=(-1)^{k-t-1}\sum_{\underline{r}\in {\frak C}_{I}^{t,k-t-1}}\chi_{Bw_{r_{k}}^{n}\cdots w_{r_{1}}^{n-k+1}B_{J_{k}}}. 
$$
Finally, substituting this into (\ref{haja}) establishes the formula. \medskip\\
3. Set ${\mathfrak C}_{\widehat{I}_{1}}^{0}=\prod_{\iota =1}^{k}\llbracket i_{\iota}-i_{1}+1,i_{\iota +1}-i_{1}\rrbracket$, and for each $t$ such that $1\leq t\leq k$,
$$
{\mathfrak C}_{\widehat{I}_{1}}^{t}\!\!\!=\!\!\left(\displaystyle\prod_{\iota =1}^{t-1}\llbracket i_{\iota\! +\!1}\!\! -\! i_{1}\! +\! 1,n\!-\!k\!+\!\iota \!+\!1\!\rrbracket\!\!\right) \!\times \! \llbracket i_{t} - i_{1}\! + 1,n\! - k + t +\! 1\rrbracket \times \!\left(\displaystyle\prod_{\iota =t+1}^{k}\!\!\!\llbracket i_{\iota}\!\! -\! i_{1}\!\! +\! 1,i_{\iota +1}\!\!-\! i_{1}\rrbracket\!\!\right) 
$$
and
$$
{\mathfrak D}_{\widehat{I}_{1}}^{t}=\prod_{\iota =1}^{t-1}\llbracket i_{\iota+1}-i_{1}+1, n-k+\iota +1\rrbracket \times \prod_{\iota =t}^{k}\llbracket i_{\iota}-i_{1}+1,i_{\iota +1}-i_{1}\rrbracket \, .
$$
We proceed as in the proof of the formula (\ref{frmle2}) above. Notice that for each $t=1,\ldots ,k$, we have ${\mathfrak C}_{\widehat{I}_{1}}^{t}={\mathfrak D}_{\widehat{I}_{1}}^{t}\amalg {\mathfrak D}_{\widehat{I}_{1}}^{t+1}$. Therefore by considering the following alternating sum :
$$
\sum_{t=1}^{k}(-1)^{k-t}\!\sum_{\underline{r}\in {\frak C}_{\widehat{I}_{1}}^{t}}\chi_{Bw_{r_{k}}^{n}\cdots w_{r_{1}}^{n-k+1}\widetilde{P}_{J_{k}}},
$$
all the sums over the ${\mathfrak D}_{\widehat{I}_{1}}^{t}$ cancel each other, except the sum over ${\mathfrak D}_{\widehat{I}_{1}}^{1}= {\mathfrak C}_{\widehat{I}_{1}}^{0}$ and the sum over ${\mathfrak D}_{\widehat{I}_{1}}^{k+1}={\mathfrak C}_{\widehat{I}_{1}}$. The alternating sum above gives then :
$$
(-1)^{k-1}\!\!\displaystyle\sum_{\underline{r}\in {\frak C}_{\widehat{I}_{1}}^{0}}\chi_{Bw_{r_{k}}^{n}\cdots w_{r_{1}}^{n-k+1}\widetilde{P}_{J_{k}}}+\sum_{\underline{r}\in {\frak C}_{\widehat{I}_{1}}}\chi_{Bw_{r_{k}}^{n}\cdots w_{r_{1}}^{n-k+1}\widetilde{P}_{J_{k}}}.
$$
Next, the two expressions above being equal, by taking the sum over ${\mathfrak C}_{\widehat{I}_{1}}^{0}$ to the other side of the equality, we get :
$$
\sum_{\underline{r}\in {\frak C}_{\widehat{I}_{1}}}\chi_{Bw_{r_{k}}^{n}\cdots w_{r_{1}}^{n-k+1}\widetilde{P}_{J_{k}}}=\sum_{t=0}^{k}(-1)^{k-t}\!\sum_{\underline{r}\in {\frak C}_{\widehat{I}_{1}}^{t}}\chi_{Bw_{r_{k}}^{n}\cdots w_{r_{1}}^{n-k+1}\widetilde{P}_{J_{k}}}.
$$
By Theorem \ref{decomp.CIPJ}, the left hand side of this equality corresponds to the decomposition of the characteristic function $\chi_{C_{\widehat{I}_{1}}\widetilde{P}_{J_{k}}}$. Therefore, by acting the element $y_{i_{1}}w_{i_{1}}$, we get :
\begin{equation}\label{2.etape1}
y_{i_{1}}w_{i_{1}}\chi_{C_{\widehat{I}_{1}}\P_{J_{k}}}=\sum_{t=0}^{k}(-1)^{k-t}\sum_{\underline{r}\in {\frak C}_{\widehat{I}_{1}}^{t}}y_{i_{1}}w_{i_{1}}\chi_{Bw_{r_{k}}^{n}\cdots w_{r_{1}}^{n-k+1}\P_{J_{k}}}.
\end{equation}
In this equality (\ref{2.etape1}), for each $t=1,\ldots ,k$, the sum over ${\frak C}_{\widehat{I}_{1}}^{t}$ is trivial in $\Sp^{k}(M)$. Indeed, when $t=1$ this follows from Proposition \ref{wsj2}.(2) and when $t$ is such that $2\leq t \leq k$ this follows from Proposition \ref{astuce}. Thus, in $\Sp^{k}(M)$, we have the equality :
\begin{equation}\label{2.etape2}
y_{i_{1}}w_{i_{1}}\chi_{C_{\widehat{I}_{1}}\P_{J_{k}}}=(-1)^{k}\sum_{\underline{r}\in {\frak C}_{\widehat{I}_{1}}^{0}}y_{i_{1}}w_{i_{1}}\chi_{Bw_{r_{k}}^{n}\cdots w_{r_{1}}^{n-k+1}\P_{J_{k}}}.
\end{equation}
Finally, by (\ref{pointe}) and since $\widetilde{T}$ is a normal subgroup of $\widetilde{N}$ (hence, for each $w\in \overline{W}$, there is $y\in \widetilde{T}\subseteq \widetilde{P}_{J_{k}}$ such that $y_{i_{1}}w=wy$), we have :
$$
y_{i_{1}}w_{i_{1}}Bw_{r_{k}}^{n}\cdots w_{r_{1}}^{n-k+1}\widetilde{P}_{J_{k}}=Bw_{i_{1}}w_{r_{k}}^{n}\cdots w_{r_{1}}^{n-k+1}\widetilde{P}_{J_{k}}.
$$
Therefore, by Lemma \ref{indicetranslation2}, for any $\underline{r}\in  {\mathfrak C}_{\widehat{I}_{1}}^{0}$ we have : 
$$
y_{i_{1}}w_{i_{1}}Bw_{r_{k}}^{n}\cdots w_{r_{1}}^{n-k+1}\widetilde{P}_{J_{k}}=Bw_{r_{k}+i}^{n}\cdots w_{r_{1}+i}^{n-k+1}\widetilde{P}_{J_{k}}.
$$ 
We conclude :
$$
\sum_{\underline{r}\in {\mathfrak C}_{\widehat{I}_{1}}^{0}}y_{i_{1}}w_{i_{1}}\chi_{Bw_{r_{k}}^{n}\cdots w_{r_{1}}^{n-k+1}\widetilde{P}_{J_{k}}}=\sum_{\underline{r}\in {\mathfrak C}_{I}^{0}} \chi_{Bw_{r_{k}}^{n}\cdots w_{r_{1}}^{n-k+1}\widetilde{P}_{J_{k}}}.
$$
This, with (\ref{2.etape2}), completes the proof. \qed

\begin{cor}\label{proprch1}
Let $I$, $i_{k+1}$ and ${\widehat I}_{1}$ as in the proposition above. Assume that $i_{k+1}=n+1$. Then, in $\Spk$ we have the identity
$$
\chi_{C_{I}\P_{J_{k}}}=(-1)^{k}y_{i_{1}}w_{i_{1}}\chi_{C_{{\widehat I}_{1}}\P_{J_{k}}}.
$$
\end{cor}
\proof Under the assumption $i_{k+1}=n+1$, in the proposition \ref{astuce2} we have ${\mathfrak C}_{I}^{k}=\emptyset$ and for each $t=1,\ldots ,k-1$
$$
{\frak C}_{I}^{t,k-t-1}\!\!=\!\!\left(\prod_{\iota=1}^{t-1}\llbracket i_{\iota}\!+\!1,n\!-\!k\!+\!\iota\! +\!1 \rrbracket \right)\!\times \!\left(\prod_{\iota=t}^{k-1}\llbracket i_{\iota +1}\!+\!1,n\!-\!k\!+\!\iota\! +\!1\rrbracket \right)\! \times \!\llbracket i_{k}\!+\!1,n\!+\!1\rrbracket, 
$$
on the other hand, we can write :
$$
\sum_{\underline{r}\in {\frak C}_{I}^{t,k-t-1}}\chi_{Bw_{r_{k}}^{n}\cdots w_{r_{1}}^{n-k+1}\P_{J_{k}}}=\sum_{(r_{1},\ldots ,r_{k-2})}\sum_{(r_{k-1},r_{k})}\chi_{Bw_{r_{k}}^{n}\cdots w_{r_{1}}^{n-k+1}\P_{J_{k}}}
$$
where $(r_{1},\ldots ,r_{k-2})$ runs through the set $(\prod_{\iota=1}^{t-1}\llbracket i_{\iota}+1,n-k+\iota +1\rrbracket )\times (\prod_{\iota=t}^{k-2}\llbracket i_{\iota +1}+1,n-k+\iota +1\rrbracket)$, and where the pair $(r_{k-1},r_{k})$ runs through  the cartesian  product of the two last intervals of ${\frak C}_{I}^{t,k-t-1}$, i.e. $\llbracket i_{k}+1,n\rrbracket \times \llbracket i_{k}+1,n+1\rrbracket$. Thus, for each $(r_{1},\ldots ,r_{k-2})$, by applying Proposition \ref{astuce} to the pair $(r_{k-1},r_{k})$ (it is clear that, if $w'=w_{r_{k-2}}^{n-2}\ldots w_{r_{1}}^{n-k+1}$, then $s_{n}w'=w's_{n}$), we deduce that 
$$
\sum_{(r_{k-1},r_{k})}\chi_{Bw_{r_{k}}^{n}\ldots w_{r_{1}}^{n-k+1}\P_{J_{k}}}\in C^{^\infty}(\G/\P_{J_{k}\cup\{n\}},\, M),
$$
then for any $t=1,\ldots ,k-1$, we have :
$$
\sum_{\underline{r}\in {\frak C}_{I}^{t,k-t-1}}\chi_{Bw_{r_{k}}^{n}\cdots w_{r_{1}}^{n-k+1}\P_{J_{k}}}\in \sum_{j=n-k+1}^{n}C^{^\infty}(\G/\P_{J_{k}\cup\{j\}},\, M).
$$ 
Consequently, from the second formula in Proposition \ref{astuce2}, it follows that in $\Sp^{k}(M)$ we have the equality :
$$
\chi_{C_{I}\P_{J_{k}}}=\sum_{\underline{r}\in {\frak C}_{I}^{0}}\chi_{Bw_{r_{k}}^{n}\cdots w_{r_{1}}^{n-k+1}\P_{J_{k}}} 
$$
which we have just proved when $k>1$ and which is obvious for $k=1$, the sets ${\frak C}_{I}^{t,k-t-1}$ being empty. By  combining this with the  third formula in Proposition \ref{astuce2}, we obtain the equality in the corollary. \qed

\section{Main theorem}

The two propositions that follow are in preparation for the proof of Theorem \ref{maintheorem}.

\begin{prop}\label{maintheorem1}
Let $k$, $1\leq k\leq n$. Let $\ph\in {\Hom}_{M}(\Sp^{k}(M),\,L)$. The map $\h_{\ph}\in {\Hom}_{M}(M[\hatI^{k}], L)$ defined by $\h_{\ph}(g(\sigma_{I},v_{0}^{o}))=\ph(g\chi_{C_{I}\P_{J_{k}}})$ for any $g\in \G$ and any $I\subseteq \De$, satisfies the harmonicity conditions.
\end{prop}
\proof Since $\G$ acts transitively on the pointed cells of a given type, we only need to show that $\h_{\ph}$ satisfies the conditions of harmonicity on the standard pointed cells. \medskip\\
{\bf(HC1)} Let $I\subseteq \De$ be such that $\De -I=\{i_{1}<\cdots <i_{k}\}$. Let ${\widehat{I}}_{1}$ be such that $\De-{\widehat{I}}_{1}=\{i_{2}-i_{1}<\cdots <i_{k}-i_{1}<n+1-i_{1}\}$. By Lemma \ref{wij} we have $(\sigma_{I},v^{o}_{i_{1}})=y_{i_{1}}w_{i_{1}}(\sigma_{\widehat{I}_{1}},v^{o}_{0})$, and by Corollary \ref{proprch1} we have :
$$
\varphi(\chi_{C_{I}\widetilde{P}_{J_{k}}})=(-1)^{k}\varphi(y_{i_{1}}w_{i_{1}}\chi_{C_{\widehat{I}_{1}}\widetilde{P}_{J_{k}}}).
$$
Therefore,
$$
{\frak h}_{\varphi}(\sigma_{I},v^{o}_{0})=(-1)^{k}{\frak h}_{\varphi}(\sigma_{I},v^{o}_{i_{1}}).
$$
Since we will need {\bf (HC3)} in order to prove {\bf (HC2)}, we will first prove {\bf (HC3)}.\medskip\\  
{\bf(HC3)} First, notice that if $(\sigma,v_{0})=(v_{0},v_{1}\ldots ,v_{k})\in \widehat{\frak I}^{k}$, since $\h_{\ph}$ satisfies {\bf (HC1)}, we have :
\begin{equation}\label{prov1}
{\frak h}_{\varphi}(\sigma, v_{0})= (-1)^{(j+1)k}{\frak h}_{\varphi}(\sigma,v_{j+1}).
\end{equation}
Notice also that for each integer $j$, $0\leq j\leq k$, we have a bijective correspondence :
\begin{equation}\label{prov53}
{\cal C}((\sigma,v_{0}), j)\; \xrightarrow{\sim} \;{\cal C}((\sigma,v_{j+1}),k),
\end{equation}
which sends the pointed cell $(\sigma', v_{0})$ (or $(\sigma',v'_{0})$ in case $j=0$) to the pointed cell $(\sigma',v_{j+1})$. Thus, by (\ref{prov1}) and (\ref{prov53}), we have : 
\begin{equation}\label{prov2}
\sum_{\sigma'\in {\cal C}((\si, v_{0}),j)}{\frak h}_{\varphi}(\sigma')=(-1)^{(j+1)k}\sum_{\sigma'\in {\cal C}((\sigma,v_{j+1}),k)}{\frak h}_{\varphi}(\sigma').
\end{equation}
Now, because of this formula (\ref{prov2}), we need only to prove that $\h_{\ph}$ satisfies {\bf (HC3)} in case $j=k$. That is, if $I\subseteq \De$ is such that $\Delta- I=\{i_{1}<\cdots <i_{k}\}$, we have to prove : 
$$
\sum_{\si\in {\cal C}((\si_{I},v_{0}^{o}),k)}{\frak h}_{\varphi}(\si)={\frak h}_{\varphi}(\si_{I},v^{o}_{0}).
$$
By the second point of Lemma \ref{prov5} and by the definition of $\h_{\ph}$ it suffices to prove that we have : 
\begin{equation}\label{prov3}
\coprod_{b\in B_{I}B_{I_{n}^{i_{k}}}/B_{I_{i_{n}}^{i_{k}}}}bC_{I_{n}^{i_{k}}}\widetilde{P}_{J_{k}}=B_{I}C_{I_{n}^{i_{k}}}\widetilde{P}_{J_{k}}=C_{I}\widetilde{P}_{J_{k}}.
\end{equation}
Both equalities are obvious (for the second equality see the definition of $C_{I}$). The union is disjoint. Indeed, take $b\in B_{I}$ such that $bC_{I_{n}^{i_{k}}}\widetilde{P}_{J_{k}}\cap C_{I_{n}^{i_{k}}}\widetilde{P}_{J_{k}}\neq \emptyset$.
By the formula (\ref{iwasawaprep2}) in the proof of Theorem \ref{decomp.CIPJ}, for each $\iota=1,\ldots k-1$, there exist $b_{\iota},b'_{\iota}\in B_{\llbracket i_{\iota}+1,n-k+\iota\rrbracket}\subseteq K^{^*}\widetilde{G}(O)$ and $p\in \widetilde{P}_{J_{k}}$ such that 
$$
bb_{k-1}\ldots b_{1}=b'_{k-1}\ldots b'_{1}p.
$$
This implies $p\in K^{^*}\widetilde{G}(O)\cap \widetilde{P}_{J_{k}}=B_{J_{k}}$ and hence $b=b'_{k-1}\ldots b'_{1}pb_{1}^{-1}\ldots b_{k-1}^{-1} \in B_{I'}$, where $I'$ is the subset of $\De$ defined as follows :
$$
I'=J_{k}\cup \bigcup_{\iota=0}^{k-1}\llbracket i_{\iota} +1, n-k+\iota\rrbracket.
$$    
But, since $n\notin I'$ then $I\cap I' \subseteq I_{n}^{i_{k}}$, we therefore have $b\in B_{I_{n}^{i_{k}}}$. \medskip\\
{\bf(HC2)} Let $I\subseteq \De$ be such that $\De -I=\{i_{1}<\cdots <i_{k}\}$. Since $\h_{\ph}$ satisfies the condition {\bf (HC3)} proved above, we can apply the lemma \ref{conditionhc3}. By this lemma, together with {\bf (HC1)}, we need only to consider the case $i_{k}=n$. Therefore we have to prove :
$$
\sum_{\si\in {\cal B}((\si_{I\cup\{n\}},v_{0}^{o}),t_{I})}\h_{\ph}(\si)=0.
$$
By the first point of Lemma \ref{prov5} and by the definition of $\h_{\ph}$ it suffices to show that in $\Sp^{k}(M)$ we have :  
\begin{equation}\label{ch21}
\sum_{b\in B_{I\cup\{n\}}/B_{I}}b\chi_{C_{I}\widetilde{P}_{J_{k}}}=0.
\end{equation}

First, notice that we have the following obvious equality, similar arguments as in the proof of (\ref{prov3}) show that the union in the right hand side is disjoint :
\begin{equation}\label{precision1}
B_{I\cup \{n\}}C_{I}\widetilde{P}_{J_{k}}=\coprod_{b\in B_{I\cup\{n\}}/B_{I}}bC_{I}\P_{J_{k}}.
\end{equation}

On the other hand, if we proceed similarly as in the proof of Theorem \ref{decomp.CIPJ} and by using Proposition \ref{generalwsj1'}, we show that we have the decomposition :  
\begin{equation}\label{ch22}
B_{I\cup\{n\}}C_{I}\widetilde{P}_{J_{k}}=\coprod_{\underline{r}\in {\frak C}_{I}^{n}}Bw_{r_{k}}^{n}w_{r_{k-1}}^{n-1}\cdots w_{r_{1}}^{n-k+1}\widetilde{P}_{J_{k}}
\end{equation}
where ${\frak C}_{I}^{n}=\left(\displaystyle\prod_{\iota =1}^{k-1}\llbracket i_{\iota}+1,n-k+\iota +1\rrbracket\right) \times \llbracket i_{k-1}+1,n+1\rrbracket$. Therefore, combining (\ref{precision1}) with (\ref{ch22}) we get :
$$
\sum_{b\in B_{I\cup \{n\}}/B_{I}}\chi_{bC_{I}\widetilde{P}_{J_{k}}}=\sum_{\underline{r}\in {\frak C}_{I}^{n}}\chi_{Bw_{r_{k}}^{n}w_{r_{k-1}}^{n-1}\cdots w_{r_{1}}^{n-k+1}\widetilde{P}_{J_{k}}}.
$$

Finally since the pair $(r_{k-1},r_{k})$ runs through the set $\llbracket i_{k-1}+1,n\rrbracket \times \llbracket i_{k-1}+1, n+1\rrbracket$ it follows from Proposition \ref{astuce} (or from Proposition \ref{wsj2} in case $k=1$) that we have :
$$
\sum_{b\in B_{I\cup \{n\}}/B_{I}}\chi_{bC_{I}\widetilde{P}_{J_{k}}} \in C^{^\infty}(\widetilde{G}/\widetilde{P}_{J_{k}\cup\{n\}},\, M). 
$$
This finishes the proof of (\ref{ch21}). \medskip\\
{\bf(HC4)} Let $I\subseteq \De$ be such that $\De-I=\{i_{1}<\cdots <i_{k}\}$. Assume $i_{k}\lneq n$ and let $i_{k+1}$ be such that $i_{k}\lneq i_{k+1}\leq n$. We have to show that :
$$
\displaystyle\sum_{j=0}^{k+1}(-1)^{j}{\mathfrak h}_{\varphi}(v^{o}_{0},v^{o}_{i_{1}},\ldots,\widehat{v^{o}_{i_{j}}},\ldots ,v^{o}_{i_{k+1}})=0.
$$
Let $\widehat{I}_{1}\subseteq \Delta$ be such that $\De - \widehat{I}_{1}=\{i_{2}-i_{1}< \cdots < i_{k+1}-i_{1}\}$, cf. Lemma \ref{wij}. By this lemma and by the definition of $\h_{\ph}$, it suffices to prove that in $\Sp^{k}(M)$ we have :
\begin{equation}\label{ch41}
y_{i_{1}}w_{i_{1}}\chi_{C_{\widehat{I}_{1}}\widetilde{P}_{J_{k}}}+\sum_{j=1}^{k+1}(-1)^{j}\chi_{C_{I_{i_{k+1}}^{i_{j}}}\widetilde{P}_{J_{k}}} =0.
\end{equation}
Recall that, for each $j=1,\ldots , k+1$, the index set which corresponds to the decomposition of the characteristic function $\chi_{C_{I_{i_{k+1}}^{i_{j}}}\widetilde{P}_{J_{k}}}$ is (cf. Theorem \ref{decomp.CIPJ}) : 
$$
{\frak C}_{I_{i_{k+1}}^{i_{j}}}=\prod_{\iota=1}^{j-1}\llbracket i_{\iota}+1,n-k+\iota+1 \rrbracket \times \prod_{\iota =j}^{k}\llbracket i_{\iota+1}+1,n-k+\iota +1 \rrbracket.
$$ 
By combining the identities (2) and (3) in Proposition \ref{astuce2}, and since we have ${\frak C}_{I_{i_{k+1}}^{i_{k}}}={\frak C}_{I}^{k}$ and ${\frak C}_{I_{i_{k+1}}^{i_{k+1}}}={\frak C}_{I}$, we obtain :
$$
y_{i_{1}}w_{i_{1}}\chi_{C_{\widehat{I}_{1}}\widetilde{P}_{J_{k}}}\!\!+(-1)^{k}\chi_{C_{I_{i_{k+1}}^{i_{k}}}\!\!\!\P_{J_{k}}}\!\!\!+(-1)^{k+1}\chi_{C_{I}\widetilde{P}_{J_{k}}}\!\!\!=\displaystyle\sum_{j=1}^{k-1}(-1)^{j}\!\!\!\!\!\!\!\!\sum_{\underline{r}\in {\mathfrak C}_{I}^{j,k-j-1}}\!\!\!\!\!\!\!\chi_{Bw_{r_{k}}^{n}\cdots w_{r_{1}}^{n-k+1}\widetilde{P}_{J_{k}}}.
$$
Substituting this into (\ref{ch41}) we conclude that :
$$
y_{i_{1}}w_{i_{1}}\chi_{C_{\widehat{I}_{1}}\widetilde{P}_{J_{k}}}+\sum_{j=1}^{k+1}(-1)^{j}\chi_{C_{I_{i_{k+1}}^{i_{j}}}\widetilde{P}_{J_{k}}}
=\displaystyle\sum_{j=1}^{k-1}(-1)^{j}\displaystyle\sum_{\underline{r}\in {\frak D}_{I_{i_{k+1}}^{i_{j}}}}\!\!\chi_{Bw_{r_{k}}^{n}\cdots w_{r_{1}}^{n-k+1}\widetilde{P}_{J_{k}}},
$$
where we have set ${\frak D}_{I_{i_{k+1}}^{i_{j}}}={\frak C}_{I_{i_{k+1}^{i_{j}}}}\amalg {\frak C}_{I}^{j,n-j-1}$. Now, it remains to show that the right hand side in the equality above is trivial in $\Sp^{k}(M)$. It is easy to see that we have :
$$
{\frak D}_{I_{i_{k+1}}^{i_{j}}}\!\!\!=\!\left(\prod_{\iota=1}^{j-1}\llbracket i_{\iota}\! +\! 1,n\! -\! k\! +\! \iota \! +\! 1\rrbracket\right)\! \times \!\left(\prod_{\iota=j}^{k-1}\llbracket i_{\iota+1}\! +\! 1,n\! -\! k\!+\!\iota\!+\! 1\rrbracket\right)\! \times \llbracket i_{k}+1,n+1\rrbracket.
$$
Notice that the cartesian product of the two last intervals in ${\frak D}_{I_{i_{k+1}}^{i_{j}}}$ is $\llbracket i_{k}+1,n\rrbracket \times \llbracket i_{k}+1,n+1\rrbracket$. Thus, we can check easily that we can apply Proposition \ref{astuce} (or Proposition \ref{wsj2} in case $k=1$) to conclude the following :
$$
\sum_{\underline{r}\in {\frak D}_{I_{i_{k+1}}^{i_{j}}}}\chi_{Bw_{r_{k}}^{n}\cdots w_{r_{1}}^{n-k+1}\widetilde{P}_{J_{k}}}\in C^{^\infty}(\widetilde{G}/\widetilde{P}_{J_{k}\cup\{n\}},\,M).
$$
This completes the proof of the assertion (\ref{ch41}). \qed

The proof of the following  proposition has been inspired by case $k=n$ which was treated by P. Schneider and J. Teitelbaum, see \cite[page 401, Lemma 10]{Schneider-T}. The computations are much more complicated for general $k$, $1\leq k\leq n$.

\begin{prop}\label{maintheorem2}
Let $\h\in {\frak Harm}^{k}(M, L)$. The map $\psi_{\h}\in {\rm Hom}_{M}(C_{c}^{^\infty}(\widetilde{G}/B_{J_{k}},\,M),\, L)$ defined by $\psi_{\h}(g\chi_{B_{J_{k}}})={\h}(g(\si_{J_{k}},v^{o}_{0}))$, for any $g\in \widetilde{G}$, vanishes on the $M[\G]$-submodule ${\frak R}_{J_{k}}$ of $C_{c}^{^\infty}(\G/B_{J_{k}},M)$.  \end{prop}
\proof 
Let $j$ be such that $n-k+1\leq j\leq n$. Since $\h$ is harmonic and then satisfies the condition {\bf (HC2)}, we conclude that :
$$
\psi_{\mathfrak h}(\chi_{B_{J_{k}\cup\{j\}}})=\sum_{b\in B_{J_{k}\cup\{j\}}/B_{J_{k}}}\psi_{\frak h}(b\chi_{B_{J_{k}}})=\sum_{b\in B_{J_{k}\cup\{j\}}/B_{J_{k}}}{\frak h}(b(\sigma_{J_{k}},v^{o}_{0}))=0.\medskip
$$
Therefore $\psi_{\h}$ is trivial on the $M[\G]$-submodule $\sum_{j=n-k+1}^{n}C_{c}^{^\infty}(\G/B_{J_{k}\cup\{j\}},M)$ of $C_{c}^{^\infty}(\G/B_{J_{k}},M)$. Let us show that $\psi_{\h}$ vanishes on the functions $\chi_{B_{J_{k}}}-\chi_{B_{J_{k}}y_{i}B_{J_{k}}}$, $0\leq i\leq n$. Since $\h$ satisfies {\bf (HC3)}, for any $I\subseteq \De$ such that $|\Delta -I|=k$, we have :
\begin{equation}\label{ch3}
\psi_{\h}(g.\chi_{C_{I}})={\h}(g(\si_{I},v^{o}_{0})).
\end{equation}
On the other hand, if $w_{i}$ is as in Lemma \ref{wij}, the assertion (\ref{pointe}) given in Remark \ref{rempointe} which follows this lemma says that $y_{i}w_{i}$ normalizes B, thus : 
\begin{equation}\label{f0}
\chi_{By_{i}B_{J_{k}}}=\chi_{By_{i}w_{i}w_{i}^{-1}B_{J_{k}}}=y_{i}w_{i}.\chi_{Bw_{i}^{-1}B_{J_{k}}}.
\end{equation}
Now, there are two cases depending on $i$ :\medskip\\
$\bullet$ $\underline{0\leq i\leq n-k+1}$ : observe that $w_{i}^{-1}$ decomposes in two factors as follows : 
$$
w_{i}^{-1}=(w_{n-i+1}^{n}\cdots w_{n-k-i+2}^{n-k+1}).(w_{n-k-i+1}^{n-k}\cdots w_{1}^{i})
$$
and that the second factor $w_{n-k-i+1}^{n-k}\cdots w_{1}^{i}$ lies in $W_{J_{k}}$ and hence in $B_{J_{k}}$. Thus, we have :
$$
Bw_{i}^{-1}B_{J_{k}}=Bw_{n-i+1}^{n}\cdots w_{n-k-i+2}^{n-k+1}B_{J_{k}}.
$$
So, if for each $\iota$, $1\leq \iota \leq k+1$, we consider $i_{\iota}=n-k-i+\iota$, and if $I\subseteq \Delta$ is such that $\Delta -I=\{i_{1} < \cdots <i_{k}\}$, we clearly have :
$$
\chi_{Bw_{i}^{-1}B_{J_{k}}}=\sum_{\underline{r}\in {\frak C}_{I}^{0}}\chi_{Bw_{r_{k}}^{n}\cdots w_{r_{1}}^{n-k+1}B_{J_{k}}},  
$$
with ${\frak C}_{I}^{0}=\prod_{\iota =1}^{k}\llbracket i_{\iota}+1,i_{\iota +1}\rrbracket$. Since $\psi_{\h}$ vanishes on $\sum_{j=n-k+1}^{n}C_{c}^{^\infty}(\G/B_{J_{k}\cup\{j\}},M)$, by the first assertion of Proposition \ref{astuce2}, we conclude that $\psi_{\frak h}$ satisfies :
\begin{equation}\label{ch2}
\psi_{\h}(\chi_{Bw_{i}^{-1}B_{J_{k}}}\!)\!=\!\psi_{\h}\!\left(\!\chi_{C_{I}}\!\!-\chi_{C_{I_{i_{k+1}}^{i_{k}}}}\!\!\!\!\!-\sum_{t=1}^{k-1}(-1)^{k-t-1}\!\!\!\!\!\!\!\!\sum_{\underline{r}\in {\frak C}_{I}^{t,k-t-1}}\!\!\!\!\!\!\!\chi_{Bw_{r_{k}}^{n}\cdots w_{r_{1}}^{n-k+1}B_{J_{k}}}\!\!\right).
\end{equation}
On the other hand, since $\h$ satisfies {\bf (HC4)}, we have :
$$
\displaystyle\sum_{j=0}^{k+1}(-1)^{j}{\h}(v^{o}_{0},v^{o}_{i_{1}},\ldots,\widehat{v^{o}_{i_{j}}},\ldots ,v^{o}_{i_{k+1}})=0.
$$
Therefore, if ${\widehat{I}}_{1}\subseteq \De$ is such that $\De-{\widehat{I}}_{1}=\{i_{2}-i_{1}<\cdots <i_{k+1}-i_{1}\}$, we have $(v_{i_{1}}^{o},\ldots ,v_{i_{k+1}}^{o})=y_{i_{1}}w_{i_{1}}(\si_{{\widehat{I}}_{1}},v_{0}^{o})$ and, by (\ref{ch3}), the identity above gives : 
\begin{equation}\label{ch4}
\psi_{\h}(y_{i_{1}}w_{i_{1}}\chi_{C_{{\widehat{I}}_{1}}})+ \sum_{t=1}^{k+1}(-1)^{t}\psi_{\h}(\chi_{C_{I_{i_{k+1}}^{i_{t}}}})=0.
\end{equation}
Combining (\ref{ch2}) with (\ref{ch4}) we conclude that :
$$
\psi_{\h}(\chi_{Bw_{i}^{-1}B_{J_{k}}})=\psi_{\frak h}\!\!\left((-1)^{k}y_{i_{1}}w_{i_{1}}\chi_{C_{{\widehat{I}}_{1}}}\!+\! \sum_{t=1}^{k-1}(-1)^{k-t}\!\!\!\!\sum_{\underline{r}\in {\frak D}_{I_{i_{k+1}}^{i_{t}}}}\!\!\chi_{Bw_{r_{k}}^{n}\cdots w_{r_{1}}^{n-k+1}B_{J_{k}}}\right),
$$
where we have set ${\frak D}_{I_{i_{k+1}}^{i_{t}}}={\frak C}_{I_{i_{k+1}}^{i_{t}}}\amalg {\frak C}_{I}^{t,k-t-1}$, and where ${\frak C}_{I_{i_{k+1}}^{i_{t}}}$ is the index set associated to the decomposition of $C_{I_{i_{k+1}}^{i_{t}}}$ (see Theorem \ref{decomp.CIPJ}). By Proposition \ref{astuce}, for each $t=1,\ldots ,k-1$, $\psi_{h}$ vanishes on the sum over ${\frak D}_{I_{i_{k+1}}^{i_{t}}}$. Therefore :
$$
\psi_{\h}(\chi_{Bw_{i}^{-1}B_{J_{k}}})=(-1)^{k}\psi_{\h}(y_{i_{1}}w_{i_{1}}\chi_{C_{{\widehat{I}}_{1}}})
$$ 
and hence, by (\ref{f0}), we have :
$$
\psi_{\h}(\chi_{By_{i}B_{J_{k}}})=(-1)^{k}\psi_{\h}(y_{i}w_{i}y_{i_{1}}w_{i_{1}}\chi_{C_{{\widehat{I}}_{1}}}).
$$
Recall that we have set $i_{1}=n-k-i+1$, therefore $y_{i}w_{i}y_{i_{1}}w_{i_{1}}C_{{\widehat I}_{1}}=y_{n-k+1}w_{n-k+1}C_{{\widehat I}_{1}}$, and hence :
\begin{equation}\label{ff1}
\psi_{\frak h}(\chi_{By_{i}B_{J_{k}}})=(-1)^{k}\psi_{\frak h}(y_{n-k+1}w_{n-k+1}\chi_{C_{{\widehat{I}}_{1}}}).
\end{equation}
Finally, since $\h$ satisfies {\bf (HC1)}, we have $\h(\si_{J_{k}},v_{n-k+1}^{o})=(-1)^{k}\h(\si_{J_{k}},v_{0}^{o})$. Thus, by Lemma \ref{wij} and by (\ref{ch3}) we get :
\begin{equation}\label{ch1}
\psi_{\h}(y_{n-k+1}w_{n-k+1}\chi_{C_{({\widehat{J}}_{k})_{1}}})=(-1)^{k}\psi_{\frak h}(\chi_{B_{J_{k}}}),
\end{equation}
where $(\widehat{J_{k}})_{1} \subseteq \De$ is such that $\De -(\widehat{J_{k}})_{1}=\{1<2<\cdots <k\}$. Note that ${\widehat{I}}_{1}=(\widehat{J_{k}})_{1}$, therefore by combining (\ref{ff1}) with (\ref{ch1}) we get $\psi_{\frak h}(\chi_{By_{i}B_{J_{k}}})=\psi_{\frak h}(\chi_{B_{J_{k}}})$. \medskip\\
$\bullet$ $\underline{n-k+2 \leq i\leq n}$ : simple calculation shows that we have $w_{i}^{-1}=ww'$ with 
$$
w=w_{n-i+1}^{2n-k-i+1}\cdots w_{2}^{n-k+2}w_{1}^{n-k+1} \quad \textrm{and} \quad w'=w_{2n-k-i+2}^{n}\cdots w_{n-k+3}^{i+1}w_{n-k+2}^{i}.
$$ 
Notice that $w'\in W_{\llbracket n-k+2,n\rrbracket}$ and that $l(w) = (n-i+1)(k-1)\;(\textrm{mod }2)$. Thus, since $\psi_{h}$ vanishes on $\sum_{j=n-k+1}^{n}C_{c}^{^\infty}(\G/B_{J_{k}\cup\{j\}},M)$, by Proposition \ref{astuce0} we have :
\begin{equation}\label{f0'}
\psi_{\frak h}(\chi_{Bw_{i}^{-1}B_{J_{k}}})=(-1)^{(n-i+1)(k-1)}\psi_{\frak h}(\chi_{Bw_{n-i+1}^{2n-k-i+1}\cdots w_{2}^{n-k+2}w_{1}^{n-k+1}B_{J_{k}}}).
\end{equation}
Next, if we consider $I\subseteq \De$ such that $\De -I=\{i_{1}<\cdots <i_{k}\}$ with the $i_{\iota}$, $1\leq \iota \leq k$,  defined as follows :
$$
i_{\iota}=\left\{\begin{array}{ll}
\iota & {\textrm si}\;\; 1\leq \iota \leq n-i+1, \\
n-k+\iota & {\textrm si}\;\; n-i+2 \leq \iota \leq n,    
\end{array}\right.
$$
then, we have the identity : 
\begin{equation}\label{f0''}
\psi_{\h}(\chi_{C_{I}})=(-1)^{n-i+1}\psi_{\h}(\chi_{Bw_{n-i+1}^{2n-k-i+1}\cdots w_{2}^{n-k+2}w_{1}^{n-k+1}B_{J_{k}}}).
\end{equation}
Indeed, since $\psi_{\h}$ vanishes on $\sum_{j=n-k+1}^{n}C_{c}^{^\infty}(\widetilde{G}/B_{J_{k}\cup\{j\}},\, M)$, by Proposition \ref{wsj2} we have : 
$$
\psi_{\frak h}(\chi_{Bw_{n-i+1}^{2n-k-i+1}\cdots w_{2}^{n-k+2}w_{1}^{n-k+1}B_{J_{k}}})=-\psi_{\frak h}\!\left(\sum_{r_{1}=2}^{n-k+2}\!\!\!\chi_{Bw_{n-i+1}^{2n-k-i+1}\cdots w_{2}^{n-k+2}w_{r_{1}}^{n-k+1}B_{J_{k}}}\!\!\!\right),
$$
and by Proposition \ref{astuce} with $(r_{1},r_{2})$ running through $\llbracket 2,n-k+2\rrbracket \times \llbracket 2,n-k+3\rrbracket$, we conclude that the right hand side in the equality above is equal to :
$$
\psi_{\frak h}\left( \sum_{r_{2}=3}^{n-k+3} \sum_{r_{1}=2}^{n-k+2}\chi_{Bw_{n-i+1}^{2n-k-i+1}\cdots w_{3}^{n-k+3}w_{r_{2}}^{n-k+2}w_{r_{1}}^{n-k+1}B_{J_{k}}}\right).
$$
Consequently, we have the identity :
$$\displaystyle
\begin{array}{c}
\psi_{\h}(\chi_{Bw_{n-i+1}^{2n-k-i+1}\cdots w_{2}^{n-k+2}w_{1}^{n-k+1}B_{J_{k}}})\\=(-1)^{2}\psi_{\h}\left( \sum_{r_{2}=3}^{n-k+3} \sum_{r_{1}=2}^{n-k+2}\chi_{Bw_{n-i+1}^{2n-k-i+1}\cdots w_{3}^{n-k+3}w_{r_{2}}^{n-k+2}w_{r_{1}}^{n-k+1}B_{J_{k}}}\right).
\end{array}
$$
By repeating this process and using Proposition \ref{astuce} successively with the pairs $(r_{2},r_{3})$, $(r_{3},r_{4})$, $\ldots$, $(r_{n-i},r_{n-i+1})$ running through the sets $\llbracket 3,n-k+3 \rrbracket \times \llbracket 3,n-k+4 \rrbracket$, $\llbracket 4,n-k+4 \rrbracket \times \llbracket 4,n-k+5 \rrbracket$, \ldots, $\llbracket n-i+1, 2n-k-i+1\rrbracket \times \llbracket n-i+1, 2n-k-i+2\rrbracket$ respectively, we get : 
$$
\begin{array}{c}
\psi_{\h}(\chi_{Bw_{n-i+1}^{2n-k-i+1}\cdots w_{2}^{n-k+2}w_{1}^{n-k+1}B_{J_{k}}})\\
=(-1)^{n-i+1}\psi_{\h}\left(\displaystyle \sum_{r_{n-i+1}=n-i+2}^{2n-k-i+2} \!\!\cdots \sum_{r_{2}=3}^{n-k+3}\sum_{r_{1}=2}^{n-k+2}\chi_{Bw_{r_{n-i+1}}^{2n-k-i+1}\cdots w_{r_{2}}^{n-k+2}w_{r_{1}}^{n-k+1}B_{J_{k}}}\!\!\right).
\end{array}
$$
Notice that the expression in the right hand side above defined by the sums over the $r_{\iota}$ is nothing else than the decomposition of $\chi_{C_{I}}$ given in Theorem \ref{decomp.CIPJ}. This proves $(\ref{f0''})$.

Finally, combining (\ref{f0'}) with (\ref{f0''}) we get $\psi_{\frak h}(\chi_{Bw_{i}^{-1}B_{J_{k}}})=(-1)^{(n-i+1)k}\psi_{\frak h}(\chi_{C_{I}})$, therefore, by using (\ref{f0}), we deduce :
\begin{equation}\label{ff1'}
\psi_{\h}(\chi_{By_{i}B_{J_{k}}})=(-1)^{(n-i+1)k}\psi_{\h}(y_{i}w_{i}\chi_{C_{I}}).
\end{equation}
On the other hand, $\h$ being harmonic, by {\bf (HC1)} we have : 
\begin{equation}\label{fff0}
\h(\si_{J_{k}},v^{o}_{i})=(-1)^{(n-i+1)k}\h(\si_{J_{k}},v^{o}_{0}).
\end{equation}
From Lemma \ref{wij}, we have $(\si_{J_{k}},v_{i}^{o})=y_{i}w_{i}(\si_{({\widehat{J_{k}}})_{i}},v_{0}^{o})$ where $({\widehat{J_{k}}})_{i}$ is nothing else than $I$. This, with (\ref{fff0}), gives :
$$
\h(y_{i}w_{i}(\si_{I},v^{o}_{0}))=(-1)^{(n-i+1)k}\h(\si_{J_{k}},v^{o}_{0}),
$$
and, by (\ref{ch3}), we deduce :
\begin{equation}\label{ch1'}
\psi_{\frak h}(y_{i}w_{i}\chi_{C_{I}})=(-1)^{(n-i+1)k}\psi_{\frak h}(\chi_{B_{J_{k}}}).
\end{equation}
Thus, combining (\ref{ff1'}) with (\ref{ch1'}) we get $\psi_{\frak h}(\chi_{By_{i}B_{J_{k}}})=\psi_{\frak h}(\chi_{B_{J_{k}}})$. \qed

\begin{theo}\label{maintheorem} For any $k$, $0\leq k\leq n$, there is an $M[\widetilde{G}]$-isomorphism 
$$
{\Hom}_{M}(\Sp^{k}(M),\,L) \cong {\frak Harm}^{k}(M,L)
$$
\end{theo}
\proof For $k=0$, by (\ref{0-special-representation}) and (\ref{0-harmonic-cochains}), both sides of the isomorphism of this theorem are canonically isomorphic to $L$. Now, let $k$, $1\leq k\leq n$. By Proposition \ref{maintheorem1}, the map which to $\ph$ associates $\h_{\ph}$ gives a well defined $M$-homomorphism 
$$
H^{k}: {\rm Hom}_{M}({\rm Sp}^{k}(M),\, L) \rightarrow {\frak Harm}^{k}(M, L).
$$
This homomorphism is clearly $\G$-equivariant. On the other hand, we have a well defined $M[\G]$-homomorphism :
$$
\begin{array}{cccc}
\Phi^{k} : & {\frak Harm}^{k}(M, L) & \longrightarrow & {\rm Hom}_{M}({\rm Sp}^{k}(M),\, L)
\end{array}
$$ 
which sends an harmonic cochain $\h$ to $\ph_{h}$ defined by $\ph_{\h}(g\chi_{B_{J_{k}}\P_{J_{k}}})=\h(g(\si_{J_{k}},v_{0}^{o}))$ for any $g\in \G$. Indeed, it is easy to check that $\Phi^{k}={\widetilde{H}_{J_{k}}}^{-1} \circ \Psi^{k}$, where the $M[\G]$-homomorphism   
$$
\begin{array}{cccc}
\Psi^{k} : & {\frak Harm}^{k}(M, L) & \longrightarrow & {\rm Hom}_{M}(C_{c}^{^\infty}(\widetilde{G}/B_{J_{k}},\,M)/{\frak R}_{J_{k}},\, L)
\end{array}
$$
which to $\h$ associates $\psi_{\h}$ is given by Proposition \ref{maintheorem2}, and $\widetilde{H}_{J_{k}}$ is the $M[\G]$-isomorphism 
$$
\begin{array}{cccc}
\widetilde{H}_{J_{k}} : & {\rm Hom}_{M}({\rm Sp}^{k}(M),\, L) & \xrightarrow{\sim} & {\rm Hom}_{M}(C_{c}^{^\infty}(\widetilde{G}/B_{J_{k}},\,M)/{\frak R}_{J_{k}},\, L),
\end{array}
$$
dual to the isomorphism $H_{J_{k}}$ given by Corollary \ref{yeap}. 

Let us prove that $\Phi^{k}$ and $H^{k}$ are inverse to each other. Let $\h\in {\frak Harm}^{k}(M,L)$ and let $\h_{\ph_{\h}}$ be its image by $H^{k} \circ \Phi^{k}$. We have ${\frak h}_{\varphi_{\frak h}}(\sigma_{J_{k}},v^{o}_{0})=\varphi_{\frak h}(\chi_{B_{J_{k}}\widetilde{P}_{J_{k}}})={\frak h}(\sigma_{J_{k}},v^{o}_{0})$. As we have $\h_{\ph_{\h}}\in {\frak Harm}^{k}(M,L)$, this proves  that $\h=\h_{\ph_{\h}}$, by the property {\bf (HC3)}.

On the other hand, if $\ph \in {\Hom}_{M}(\Sp^{k}(M),\, L)$ and if $\ph_{\h_{\ph}}$ is its image by $\Phi^{k} \circ H^{k}$, then we have $\varphi_{{\frak h}_{\varphi}}(\chi_{B_{J_{k}}\widetilde{P}_{J_{k}}})={\frak h}_{\varphi}(\sigma_{J_{k}},v^{o}_{0})=\varphi(\chi_{B_{J_{k}}\widetilde{P}_{J_{k}}})$. 
We are done. \qed

LMNO-UMR, Universit\'e de Caen, Campus II, boulevard Mar\'echal Juin, BP 5186, 14032 Caen cedex, France.\\ 
e-mail adress : amrane@math.unicaen.fr (ou : amrane@picard.ups-tlse.fr).

\end{document}